\documentclass[aop,citesort,MSNbibl,dvips]{arximspdf}

\doi{10.1214/10-AOP563}
\volume{39}
\issue{5}
\pubyear{2011}
\firstpage{1668}
\lastpage{1701}

\makeatletter

\newtheorem{theorem}{Theorem}[section]
\newtheorem{conj}[theorem]{Conjecture}
\newproclaim{defn}[theorem]{Definition}
\newproclaim{problem}[theorem]{Problem}
\newproclaim{example}[theorem]{Example}
\newtheorem{prop}[theorem]{Proposition}

\newcommand{\E}{{\mathbb E}}
\newcommand{\RR}{{\mathbb R}}
\newcommand{\Z}{{\mathbb Z}}
\newcommand{\T}{{\mathbb T}}
\newcommand{\R}{{\mathbb R}}
\newcommand{\C}{{\mathbb C}}
\newcommand{\Hb}{{\mathbb H}}
\newcommand{\Pb}{{\mathbb P}}

\newcommand{\eps}{\varepsilon}
\newcommand{\dist}{\operatorname{dist}}
\newcommand{\Aut}{\operatorname{Aut}}
\newcommand{\Stab}{\operatorname{Stab}}
\newcommand{\cA}{\mathcal{A}}
\newcommand{\cB}{\mathcal{B}}
\newcommand{\FSF}{\mathit{FSF}}
\newcommand{\WSF}{\mathit{WSF}}
\newcommand{\FMSF}{\mathit{FMSF}}
\newcommand{\WMSF}{\mathit{WMSF}}

\makeatother

\begin{document}
\begin{frontmatter}

\title{Percolation beyond $\Z^{\bolds{d}}$: The contributions of~Oded~Schramm}
\runtitle{Percolation beyond $\Z^d$}

\begin{aug}
\author[A]{\fnms{Olle} \snm{H\"{a}ggstr\"{o}m}\corref{}\thanksref{t1}\ead[label=e1]{olleh@chalmers.se}%
\ead[label=u1,url]{http://www.math.chalmers.se/\textasciitilde olleh/}}

\runauthor{O. H\"{a}ggstr\"{o}m}

\affiliation{Chalmers University of Technology}

\dedicated{This paper is dedicated to the memory of Oded Schramm}

\address[A]{Department of Mathematical Sciences \\
Chalmers University of Technology \\
412 96 G\"oteborg \\
Sweden \\
\printead{e1}\\
\printead{u1}} 
\end{aug}

\thankstext{t1}{Supported by the G\"oran Gustafsson Foundation for Research
in the Natural Sciences and Medicine.}

\received{\smonth{5} \syear{2010}}

%
\begin{abstract}
Oded Schramm (1961--2008) influenced greatly the development
of percolation theory beyond the usual $\Z^d$ setting; in particular, the
case of nonamenable lattices. Here, we review some of his work in this field.
\end{abstract}

%

\end{frontmatter}

\section{Introduction}
\label{sect:intro}

Oded Schramm was born in 1961 and died in a hiking accident in 2008, in
what otherwise seemed to be the middle of an extraordinary mathematical
career. Although he made seminal contributions to many areas of
mathematics in general and probability in particular, I will here
restrict attention to his work on percolation processes taking place on
graph structures more exotic than the usual $\Z^d$ setting. The title I
have chosen alludes to the short but highly influential paper
\textit{Percolation beyond $\Z^d$}, \textit{many questions and a few answers}
from 1996 by Itai Benjamini and Oded Schramm \cite{BS96}.

I need to point out, however, that there are at least two respects in
which I will fail to deliver on what my chosen title suggests. First,
I will not come anywhere near an exhaustive exposition of Oded's
contributions to the field. All I can offer is a personal and highly
subjective selection of highlights. Second, Oded was a very
collaborative mathematician, and I will make no attempt (if it even
makes sense) at identifying his individual contributions as opposed to
his coauthors'. Suffice it to say that everyone who worked with him
knew him as a very generous person and as someone who would not put his
name on a paper unless he had contributed at least his fair share. I
will just quote one recollection from Oded's long-time collaborator and
friend Russ Lyons:

\begin{quote}
To me, Oded's most distinctive mathematical talent was his
extraordinary clarity of thought, which led to dazzling proofs and
results. Technical difficulties did not obscure his vision. Indeed,
they often melted away under his gaze. At one point when the four of
us [Oded, Russ, Itai Benjamini and Yuval Peres] were working on uniform
spanning forests, Oded came up with a brilliant new application of the
Mass-Transport Principle. We were not sure it was kosher, and I still
recall Yuval asking me if I believed it, saying that it seemed to be
``smoke and mirrors.'' However, when Oded explained it again, the smoke
vanished \cite{MOSblog}.
\end{quote}

\noindent Following the spirit of the aformentioned paper \cite{BS96}, I will take
``percolation beyond $\Z^d$'' to mean percolation process
that are not naturally
thought of as embedded in $d$-dimensional Euclidean space. This excludes
contributions by Oded not only to the theory of percolation on $\Z^d$
(such as \cite{BHS}) but also to percolation on the triangular lattice
(such as \cite{S01} and \cite{SS})
and to continuum percolation in $\RR^d$ (such as \cite{BS98}).

Other parts of Oded's work are discussed in the papers by Angel, Garban
and Rohde in the present volume. For further reactions to Oded's untimely
death, and memories of his life and work, see, for instance,
Lyons~\cite{L08}, H\"aggstr\"om \cite{H08} and Werner \cite{W}, as
well as
the blog \cite{MOSblog}.

\section{How it began} \label{sect:BS-beyond}

It will be assumed throughout that
$G=(V,E)$ is an infinite but locally finite connected graph. In
\textit{i.i.d. site percolation} on $G$ with retention parameter $p \in[0,1]$,
each vertex $v \in V$ is declared open
(retained, value $1$) with probability $p$
and closed (deleted, value $0$) with the remaining probability $1-p$,
and this is done independently for different vertices. Alternatively,
one may consider \textit{i.i.d. bond percolation}, which is similar except
that it is the edges rather than the vertices that are declared open or closed.
Write $\Pb_{p,\mathrm{site}}$ and $\Pb_{p,\mathrm{bond}}$ for the resulting probability
measures on $\{0,1\}^V$ and $\{0,1\}^E$, respectively.
The choice whether to study bond or site percolation is often (but not always)
of little importance and largely a matter of taste. In either case, focus
is on the connectivity structure of the resulting random subgraph of $G$.
Of particular interest is the possible occurrence of an infinite connected
component---an infinite cluster, for short. The probability under
$\Pb_{p,\mathrm{site}}$ or $\Pb_{p,\mathrm{bond}}$ of having an infinite cluster is
always $0$
or~$1$, and increasing in $p$. This motivates defining the site
percolation critical value
\[
p_{c,\mathrm{site}}(G) = \inf\{p\dvtx\Pb_{p, \mathrm{site}}
( \exists\mbox{ an infinite cluster})>0\}
\]
and the bond percolation critical value $p_{c,\mathrm{bond}}(G)$ analogously.

By far the most studied case is where $G=(V,E)$ is the $\Z^d$ lattice
with $d \geq2$, meaning
that $V= \Z^d$ and $E$ consists of all pairs of Euclidean nearest neighbors.
Some selected landmarks in the history of
percolation are the 1960\vadjust{\eject} result of Harris
\cite{Har} that $p_{c,\mathrm{bond}}(\Z^2) \geq\frac{1}{2}$; the 1980 result
of Kesten~\cite{K} that $p_{c,\mathrm{bond}}(\Z^2) = \frac{1}{2}$; the 1987 result
of Aizenman, Kesten and Newman~\cite{AKN} establishing
uniqueness of the infinite cluster for arbitrary $d$;
and the strikingly short and beautiful alternative proof from
1989 by Burton and Keane~\cite{BK} of the same result.
See Grimmett~\cite{Gr} and Bollob\'as and Riordan~\cite{BR} for
introductions to percolation theory with
emphasis on the $\Z^d$ case.

Benjamini and Schramm \cite{BS96} were of course not the first to
study percolation on more exotic graphs and lattices. The case where
$G$ is the $(d+1)$-regular tree $\T_d$ had been well understood for a long
time, essentially because it can be seen as a Galton--Watson process.
Lyons \cite{L90,L92} had studied percolation on general trees, and
Grimmett and Newman \cite{GN} had considered percolation on the Cartesian
product $\T_d \times\Z$ [the Cartesian product $G=(V,E)$ of
two graphs $G_1=(V_1, E_1)$ and $G_2=(V_2, E_2)$ has vertex set
$V=V_1 \times V_2$, and an edge connecting $(x_1, x_2)$ and $(y_1, y_2)$
iff $x_1=y_2$ are identical and $x_2$ and $y_2$ are neighbors or vice versa].
However, it was only with the publication of Benjamini and Schramm
\cite{BS96}
that a systematic study of percolation beyond $\Z^d$ began to take off
toward anything like escape velocity. Clearly, they managed to
find just the right time for
launching the kind of informal research programme that
their paper proposes. It should be noted, however, that a large part of
the meaning of ``just the right time'' in this context is simply
``soon after Oded Schramm had been drawn into probability
theory.''

When, as in \cite{BS96}, we focus on the occurrence and properties of
infinite clusters for i.i.d. site (or bond) percolation on a graph $G=(V,E)$,
a first basic issue is of course whether such clusters occur at
all for any nontrivial value of $p$, that is, whether $p_{c,\mathrm{site}}(G)<1$
[or $p_{c,\mathrm{bond}}(G)<1$]. Benjamini and Schramm conjecture
the following, where, for $v \in V$, $B(v,n)$ denotes the set of vertices
$w \in V$ such that $\dist_G(v,w) \leq n$, and $\dist_G$ is graph-theoretic
distance in $G$.
\begin{conj}[(Conjecture~2 in \cite{BS96})] \label{conj:nontrivial}
If $G$ is a quasi-transitive graph such
that for some (hence any) $v \in V$, $|B(v,n)|$ grows faster
than linearly, then $p_{c,\mathrm{site}}(G)<1$.
\end{conj}

Here, of course, we need to define quasi-transitivity of a graph (the
term used in \cite{BS96} was \textit{almost transitive}, but the mathematical
community quickly decided that quasi-transitive was preferable).
\begin{defn}
Let $G=(V,E)$ be an infinite locally finite connected graph. A bijective
map $f\dvtx V \rightarrow V$ such that $\langle f(u), f(v) \rangle\in E$
if and only if $\langle u,v \rangle\in E$ is called a
\textit{graph automorphism} for $G$. The graph $G$ is said to
be \textit{transitive} if for any $u,v \in V$ there exists a graph
automorphism $f$ such that $f(u)=v$. More generally, $G$ is
said to be \textit{quasi-transitive} if there is a $k<\infty$ and a partitioning
of $V$ into $k$ sets $V_1, \ldots, V_k$ such that for $i=1, \ldots,
k$ and
any $u,v\in V_i$ there exists a graph
automorphism $f$ such that $f(u)=v$.
\end{defn}

An important subclass of transitive graphs is the class of graphs
arising as the Cayley graph of a finitely generated group.
It may be noted that for quasi-transitive graphs (and more generally
for bounded degree graphs; cf.~\cite{H00}) we have $p_{c,\mathrm{bond}}<1$ iff
$p_{c,\mathrm{site}}<1$, so Conjecture \ref{conj:nontrivial} may equivalently be
phrased for bond percolation. Benjamini and Schramm found a short and
elegant proof of the conjecture for the special case of so-called
nonamenable graphs:\vspace*{-2pt}
\begin{defn} \label{defn:Cheeger}
The \textit{isoperimetric constant} $h(G)$ of a graph $G=(V,E)$ is
defined as
\[
h(G) = \inf_S \frac{|\partial S|}{|S|} ,
\]
where the infimum ranges over all finite nonempty subsets of $V$, and
$\partial S= \{u \in V\setminus S\dvtx\exists v \in S \mbox{ such that }
\langle u,v \rangle\in E\}$. The graph $G$ is said to be \textit{amenable}
if $h(G)=0$; otherwise, it is said to be \textit{nonamenable}.\vspace*{-2pt}
\end{defn}

(Sometimes, as in \cite{BS96},
$h(G)$ is also called the Cheeger constant.)\vspace*{-2pt}
\begin{theorem}[(Theorem 2 in \cite{BS96})]
Any nonamenable graph $G$ satisfies $p_{c,\mathrm{site}}(G)<1$. In fact,
%
%
\begin{equation} \label{eq:nontrivial_pc}
p_{c,\mathrm{site}}(G) \leq\frac{1}{h(G)+1}.
\end{equation}
\vspace*{-2pt}
\end{theorem}
\begin{pf}
Fix $p$ and
a vertex $\rho\in V$, and consider the following sequential procedure
for searching the open cluster containing $\rho$. If $\rho$ is open,
set $S_1=\{\rho\}$, otherwise stop. At each integer time $n$, check the
status (open or closed) of some thitherto unchecked
vertex $v$ in $\partial S_{n-1}$; if $v$ is open we set $S_n=S_{n-1}
\cup\{v\}$,
if $v$ is closed we set $S_n=S_{n-1}$, while if no such $v$ can be found
the procedure terminates. Define $X_0=0$ and $X_n=|S_n|$ for $n \geq1$,
and note that $\{X_0, X_1, \ldots\}$ is a random walk whose i.i.d. increments
take value $0$ with probability $1-p$ and $1$ with probability $p$, stopped
at some random time. If $G$ has isoperimetric constant $h(G)$, then the
random walk can stop only when $\frac{n-X_n}{X_n} \geq h(G)$, that is, when
%
%
\begin{equation} \label{eq:SLLN_argument}
\frac{X_n}{n} \leq\frac{1}{h(G)+1} .
\end{equation}
But the random walk has drift $p$,
so when $p>\frac{1}{h(G)+1}$ the Strong Law of Large Numbers
implies that with positive probability
$\frac{X_n}{n}$ never satisfies (\ref{eq:SLLN_argument}), in which case
the walk never stops and $\rho$ belongs to an infinite cluster.\vspace*{-2pt}
\end{pf}

It may be noted that the result is sharp in the sense that
the bound (\ref{eq:nontrivial_pc}) holds with equality for the
tree $\T_d$, for which $p_{c,\mathrm{site}}(\T_d)=\frac{1}{d}$ and
$h(\T_d)=d-1$.

Once $p_{c,\mathrm{site}}(G)<1$, we know that i.i.d. site percolation on $G$
has two distinct phases: for $p<p_{c,\mathrm{site}}$ there is no infinite cluster, while
for $p>p_{c,\mathrm{site}}$ there is. But what happens \textit{at} the critical value?
This has long been a~central issue in percolation theory. When $G$
is the $\Z^d$ lattice with $d\geq2$, the consensus belief among
percolation theorists is that there is no infinite cluster at criticality; this is known
for $d=2$ (Russo~\cite{Rus}) and $d\geq19$ (Hara and Slade~\cite{HS}),
but the general case remains open. Benjamini and Schramm suggest that
``it might be beneficial to study the problem in other
settings.''\vspace*{2pt}\looseness=1
\begin{conj}[(Conjecture 4 in \cite{BS96})] \label{conj:no_perc_at_criticality}
$\!\!$For any quasi-transitive graph~$G$ with $p_{c,\mathrm{site}}(G)<1$, there is
a.s. no infinite cluster at criticality.
\end{conj}

This remains open (as of course it must be as long as the $\Z^d$ case
with $3 \leq d \leq18$ stays unsolved), but Benjamini and Schramm
were soon to be involved in remarkable progress toward proving it;
see Section \ref{sect:at_criticality}. The quasi-transitivity
condition cannot be dropped, as it is easy to construct graphs with
a nontrivial critical value which nevertheless percolate at criticality
(a~tree growing slightly faster than a binary tree will do;
see, e.g.,~\cite{L92}).

Another natural next question,
once the existence of infinite clusters at nontrivial values of $p$
(i.e., $p_{c,\mathrm{site}}<1$)
is established, concerns how many infinite clusters there can be.
Benjamini and Schramm \cite{BS96} noted that the argument of Newman
and Schulman \cite{NS} for showing that the number of infinite clusters
is, for fixed $p$, an a.s. constant which must equal $0$, $1$ or
$\infty$ extends to the setting of quasi-transitive graphs.
They furthermore saw that the argument of Burton and Keane \cite{BK}
for ruling out infinitely many infinite clusters extends to
amenable quasi-transitive graphs (a similar observation was made earlier
in Theorem 1$'$ of
Gandolfi, Keane and Newman \cite{GKN}). In other words, we get the following.\vspace*{2pt}
\begin{theorem}
For any amenable quasi-transitive graph $G$
and any $p \in[0,1]$, the number
of infinite clusters produced by i.i.d. site or bond
percolation on $G$ with parameter $p$ is either $0$ or $1$ a.s.
\end{theorem}

In contrast, i.i.d. site percolation
the regular tree $\T_d$ with $d\geq2$ exhibits infinitely many
infinite clusters for all $p \in(p_{c,\mathrm{site}}, 1)$. Also, the
$\T_d \times\Z$ example studied by Grimmett and Newman \cite{GN}
exhibits the same phenomenon when~$p$ is above but sufficiently
close to $p_{c,\mathrm{site}}$ (Grimmett and Newman showed this for large~$d$, and
later Schonmann \cite{Scho99b} indicated how to do it for all $d \geq2$).
Benjamini and Schramm \cite{BS96} conjectured
that the Burton--Keane argument is sharp in the sense that
whenever $G$ is quasi-transitive and nonamenable, uniqueness of
the infinite cluster fails for $p$ above but sufficiently close to
$p_{c,\mathrm{site}}$. In terms of the so-called uniqueness critical value
\begin{eqnarray*}
p_{u,\mathrm{site}} &=& p_{u,\mathrm{site}}(G)\\
&=&
\inf\{p \in[0,1]\dvtx\mbox{i.i.d. site percolation
on $G$ with parameter $p$} \\
&&\hspace*{92.67pt} \mbox{ produces a.s. a unique infinite cluster}\} ,
\end{eqnarray*}
their conjecture reads as follows.
\begin{conj}[(Conjecture 6 in \cite{BS96})] \label{conj:nonuniqueness_phase}
For any nonamenable quasi-transitive graph $G$, we have
$p_{c,\mathrm{site}}(G) < p_{u,\mathrm{site}}(G)$.
\end{conj}

In the spirit of the Grimmett--Newman result mentioned above,
they proved that for any quasi-transitive graph $G$, the product
graph $\T_d \times G$ satisfies $p_{c,\mathrm{site}}(G) < p_{u,\mathrm{site}}(G)$ provided
$d$ is large enough (this is Corollary 1 in \cite{BS96}).

Conjecture \ref{conj:nonuniqueness_phase} has stimulated further research
as well. For instance, Pak and Smirnova-Nagnibeda
\cite{PSN} proved in the case of bond percolation
that it holds for the Cayley graph of any nonamenable
group provided an appropriate choice of generators. Lalley \cite{Lal}
and a later paper by Benjamini and Schramm \cite{BS01} proved
it for certain classes of nonamenable planar graphs; this will be discussed
in more detail in Section \ref{sect:planar}.

By definition of $p_{u,\mathrm{site}}$ and the Newman--Schulman $0$--$1$--$\infty
$ law, we have infinitely many infinite clusters a.s. for any
$p \in(p_{c,\mathrm{site}}, p_{u,\mathrm{site}})$.
It would be nice to add that uniqueness holds for all $p \in(p_{u,\mathrm{site}},1)$,
but for this we need the monotonicity property that if $p_1<p_2$ and
uniqueness holds a.s. at parameter value $p_1$, then it
holds at $p_2$ as well; this is part of Question 5 in \cite{BS96}.
The required monotonicity was proved by
H\"aggstr\"om and Peres \cite{HP} for quasi-transitive
graphs under the additional assumption of unimodularity (see Definition
\ref{defn:unimodularity} below), and by Schonmann \cite{Scho} without this
additional assumption.

The other part of Question 5 in \cite{BS96} concerns, in the case where
$G$ is quasi-transitive with
$p_{c,\mathrm{site}}<p_{u,\mathrm{site}}<1$, the number of infinite cluster at $p=p_{u,\mathrm{site}}$:
one or infinitely many? Somewhat surprisingly, the answer
turned out to depend on the choice of $G$. For the Grimmett--Newman
example, Schonmann \cite{Scho99b} showed that there are infinitely
many infinite clusters at the uniqueness critical point $p_{u,\mathrm{site}}$, a result
that Peres \cite{Per} extended to more general product graphs. In contrast,
Benjamini and Schramm \cite{BS01} showed that for planar nonamenable graphs
with one end, there is a unique infinite cluster at $p_{u,\mathrm{site}}$; see
Section \ref{sect:planar} again.

There is a good deal more to say about the \textit{Percolation beyond
$\Z ^d$} paper~\cite{BS96}, but I must move on to some of Oded Schramm's later contributions.
The paper's influence will be evident from the coming sections,
but see also Benjamini and Schramm \cite{BSprogress}, Lyons \cite{L00}
and H\"aggstr\"om and Jonasson \cite{HJ06} for partially overlapping surveys
of what happened in the wake of the paper.

\section{Invariant percolation and mass transport}
\label{sect:BLPS-gip}

Soon after finishing the \textit{Percolation beyond $\Z^d$} paper
\cite{BS96}, Itai Benjamini and Oded Schramm joined forces with
Russ Lyons and Yuval Peres (this quartet of authors will
appear frequently in what follows, and will be abbreviated BLPS). In
\cite{BLPSgip}, they broadened the scope compared
to \cite{BS96} by considering percolation processes on
quasi-transitive graphs in a more general situation than the i.i.d.,
namely \textit{automorphism invariance}.
\begin{defn}
Let $G=(V,E)$ be a quasi-transitive graph and let $\Aut(G)$ denote the
group of graph automorphisms of $G$. A $\{0,1\}^V$-valued random
object $X$ is called a site percolation for $G$, and it is said to
be \textit{automorphism invariant} if for any $n$, any $v_1, \ldots, v_n
\in V$,
any $b_1, \ldots, b_n \in\{0,1\}$ and any $\gamma\in\Aut(G)$, we
have
\[
\Pb\bigl(X(\gamma v_1)=b_1, \ldots, X(\gamma v_n)=b_n\bigr)=
\Pb\bigl(X(v_1)=b_1, \ldots, X(v_n)=b_n\bigr) .
\]
\end{defn}

Automorphism invariance of a bond percolation for $G$ is defined
analogously.
In fact, much of the work in \cite{BLPSgip} concerns an even more
general setting, namely invariance under certain kinds of subgroups
of $\Aut(G)$. Here, for simplicity, I~will restrict attention
to the case of invariance under the full automorphism group $\Aut(G)$.

There are plenty of automorphism invariant percolation processes beyond
i.i.d. that arise naturally. Examples in the site precolation case
include certain Gibbs distributions
for spin systems such as the Ising model, and certain equilibrium measures
for interacting particle systems such as the voter model. In the bond
percolation case, they include the random-cluster model~\cite{Gr06} as
well as the random spanning forest models to be discussed in
Section~\ref{sect:forests}.

Amongst the most important contributions of BLPS \cite{BLPSgip} is the
introduction of the so-called \textit{Mass-Transport Principle} in
percolation theory, and the beginning of a systematic exploitation
of it for understanding the behavior of percolation processes.
(This was partly inspired by an application in H\"aggstr\"om \cite{H97}
of similar ideas in the special case where $G$ is a regular tree.)
As a kind of warm-up for readers unfamiliar with the mass-transport
technique, let me suggest a very simple toy problem.

\begin{problem} \label{prob:minimal_clusters}
Given a transitive graph $G=(V,E)$,
does there exist an automorphism invariant bond percolation
process which produces, with positive probability, some infinite
open cluster consisting of a single self-avoiding path which is
infinite in just one direction? (We call such a self-avoiding
path \textit{uni-infinite}.) In other words, this open cluster should
consist of a single vertex of degree $1$ in the cluster,
while all the other (infinitely
many) vertices of the cluster should have degree $2$.
\end{problem}

Call an infinite cluster \textit{slim} if it is of the desired kind.
(In a sense, a~slim infinite cluster is the smallest
infinite cluster there can be.) Also,
given an automorphism invariant bond percolation process $X$ taking
values in $\{0,1\}^E$, we define random variables $\{Y(v)\}_{v\in V}$ as
follows. If $v$ does not belong to a slim infinite cluster
in~$v$, we set $Y(v)=0$; otherwise, we let $Y(v)$ be one plus the distance
in the slim cluster from $v$ to the one endpoint of this cluster. For
$k=0,1, \ldots,$ write $\alpha(v,k)= \Pb(Y(v)=k)$. Automorphism invariance
ensures that this is independent of the choice of $v$, so we may write
$\alpha(k)$ for $\alpha(v,k)$.

When $G=(V,E)$ is the $\Z^d$ lattice, we can argue as follows. Write
$\Lambda_n$
for the box $\{-n, -n+1, \ldots, n\}^d \subset V$. The expected number of
vertices $v \in\Lambda_n$ with $Y(v)=1$ is $(2n+1)^d\alpha(1)$. But
for any $k$, and any vertex $v$ with $Y(v)=1$, there must be a corresponding
vertex $u$ with $Y(u)=k$ within distance $k-1$ from $v$. Hence,
the expected number of vertices $v \in\Lambda_{n+k-1}$ with $Y(v)=k$ is
at least $(2n+1)^d\alpha(1)$, so
\[
\bigl(2(n+k)-1\bigr)^d \alpha(k) \geq(2n+1)^d\alpha(1) ,
\]
and sending $n \rightarrow\infty$ yields $\alpha(k) \geq\alpha(1)$.
Similarly, $\alpha(1) \geq\alpha(k)$, so that in fact $\alpha
(1)=\alpha(k)$,
and since $k$ was arbitrary we have $\alpha(1)=\alpha(2)=\alpha
(3)=\cdots.$
But $\sum_{k=1}^\infty\alpha(k) \leq1$, so $\alpha(k)$ must be $0$ for
each $k$, whence slim infinite clusters do not occur in the $G=\Z^d$ case.

The crucial property of the $\Z^d$ lattice that makes the argument
work is that $\frac{|\Lambda_{n+k}|}{|\Lambda_n|} \rightarrow1$ as
$n \rightarrow\infty$. Hence, the argument is easily extended to the more
general case where $G$ is transitive and amenable.

But what about the case where $G$ is nonamenable? Now the argument does not
generalize, and in fact the following example gives us problems. And
\textit{end} in a graph $G=(V,E)$ is an equivalence class of uni-infinite
self-avoiding paths in $X$, with two paths equivalent if for all
finite $W \subset V$ the paths are eventually in the same connected
component of the graph obtained from $G$ by deleting all $v\in W$.
\begin{example}[(Trofimov's graph \cite{T})] \label{ex:Trofimov}
Consider the regular binary\break tree~$\T_2$,~and fix an end $\xi$ in
this tree. For each vertex $v$ in the tree, there is a unique uni-infinite
self-avoiding path from $v$ that belongs to $\xi$. Call the first vertex
after $v$ on this path the $\xi$-\textit{parent} of $v$, and call the
other two neighbors of $v$ its $\xi$-\textit{children}. The
$\xi$-\textit{grandparent} of $v$ is defined similarly in the obvious way.
Let $G=(V,E)$ be the graph that arises by taking
$\T_2$ and adding, for each vertex $v$, an extra edge connecting $v$
to its $\xi$-grandparent.
\end{example}

Clearly, Trofimov's graph $G$ is transitive, and it also inherits
the nonamenability property of $\T_2$.
It turns out that on $G$, it is possible to construct an automorphism
invariant bond percolation exhibiting slim infinite
clusters:

\begin{example} \label{example:Trofimov_percolation}
Let $G=(V,E)$ be Trofimov's graph, and consider the following
automorphism invariant bond percolation on $G$: each $v\in V$ will have
an open edge to exactly one of its $\xi$-children, and for each $v$
independently, toss a fair coin to decide which $\xi$-child to connect to.
All grandparent--grandchild edges are closed.
(To see that this bond percolation is indeed automorphism invariant, it
is necessary, but easy, to check that the end $\xi$ can be identified by
just looking at the graph structure of $G$.) This produces a percolation
configuration in which a.s. each $v$ sits in a slim infinite cluster.
From $v$ the open path extends downward (i.e., away from $\xi$) infinitely,
and upward a geometric$(\frac{1}{2})$ number of steps.
\end{example}

So perhaps slim infinite clusters can arise as soon as $G$ is nonamenable?
In fact, no. It turns out that
the mass-transport method of BLPS \cite{BLPSgip} applies to rule out
slim infinite clusters also in the nonamenable case, as long
as the graphs satisfy the additional assumption of \textit{unimodularity}.
Unimodularity holds for all specific examples considered so far except
for Trofimov's graph. It holds for Cayley graphs in general, and I daresay
it tends to hold for most transitive graphs that are not constructed
for the explicit purpose of being nonunimodular.
\begin{defn} \label{defn:unimodularity}
Let $G=(V,E)$ be a quasi-transitive graph with automorphism group
$\Aut(G)$. For $v \in V$, the \textit{stabilizer}
of $v$ is defined as $\Stab(v)= \{\gamma\in\Aut(G)\dvtx\gamma v = v\}$.
The graph $G$ is said to be \textit{unimodular} if for all $u,v\in V$
in the
same orbit of $\Aut(G)$ we
have the symmetry
\[
|{\Stab(u)v}|=|{\Stab(v)u}|.
\]
\end{defn}

[Note how Trofimov's graph fails to be unimodular: for two vertices $u$
and~$v$ such that $u$ is the $\xi$-parent of $v$, we get $|{\Stab(u)v}|\,{=}\,2$
but \mbox{$|{\Stab(v)u}|\,{=}\,1$}. Each vertex has two children but just one
parent.]
\begin{prop} \label{prop:no_minimal_infinite_clusters}
In automorphism invariant bond percolation on a~qua\-si-transitive
unimodular graph $G$, there is a.s. no slim infinite cluster.
\end{prop}

This, we will find, is an easy consequence of the Mass-Transport Principle
of BLPS \cite{BLPSgip}. For
an automorphism invariant site (or bond) percolation
on a quasi-transitive graph $G=(V,E)$, let $\mu$ be the corresponding
probability measure on $\{0,1\}^V$ (or on $\{0,1\}^E$).
Consider a nonnegative func\-tion~$m(u,v,\omega)$ of three variables:
two vertices $u,v \in V$ and the percolation configuration $\omega$ taking
values in $\Omega= \{0,1\}^V$ (or $\Omega= \{0,1\}^E$).
Intuitively, we should think
of $m(u,v,\omega)$ as the mass transported from $u$ to $v$ given
the configuration $\omega$. We assume that $m(u,v,\omega)=0$ unless
$u$ and $v$ are in the same orbit of $\Aut(G)$, and furthermore that
$m(\cdot, \cdot, \cdot)$ is
invariant under the diagonal action of $\Aut(G)$, meaning that
$m(u,v,\omega)= m(\gamma u, \gamma v, \gamma\omega)$ for all~$u,v,
\omega$ and $\gamma\in\Aut(G)$.
\begin{theorem}[(The Mass-Transport Principle, Section 3 in \cite
{BLPSgip})] \label{thm:mass-transport}
$\!\!$Given~$G$, $\mu$ and $m(\cdot, \cdot, \cdot)$ as above, let
\[
M(u,v) = \int_\Omega m(u,v,\omega) \,d \mu(\omega)
\]
for any $u,v\in V$. If $G$ is unimodular,
then the expected total mass transported out of
any vertex $v$ equals the expected mass transported into $v$, that is,
%
%
\begin{equation} \label{eq:mass_transport_equilbrium}
\sum_{u \in V} M(v,u) = \sum_{u \in V} M(u,v) .
\end{equation}
\end{theorem}

The Mass-Transport Principle as stated here fails if $G$ is not unimodular.
[To see this for Trofimov's graph, we can consider the the mass transport
in which each vertex simply sends unit mass to its $\xi$-parent, regardless
of the percolation configuration. Then each vertex sends mass $1$ but
receives mass $2$, thus violating~(\ref{eq:mass_transport_equilbrium}).]
In fact, BLPS \cite{BLPSgip} did state a version of the
Mass-Transport Principle that holds also in the nonunimodular case;
this involves a reweighting of the mass sent from $u$ to $v$ by a factor
that depends on $\frac{|{\Stab(u)v}|}{|{\Stab(v)u}|}$. But it is in
the unimodular case that the Mass-Transport Principle has turned out
most useful, and for simplicity we stick to this case.

The proof of the Mass-Transport Principle is particularly simple in
the case where $G$ is the Cayley graph of a finitely generated group
$H$, so
here I will settle for that case only:\vspace*{-2pt}
\begin{pf*}{Proof of Theorem \ref{thm:mass-transport} in the Cayley
graph case}
For $u,v \in V$, we also have that $u$ and $v$ are group elements of $H$,
and that there is a~unique element $h=uv^{-1}\in H$ such that
$u=hv$. This gives
\begin{eqnarray*}
\sum_{u \in V} M(v,u) & = & \sum_{h \in H} M(v, hv)
= \sum_{h \in H} M(h^{-1}v, v) \\
& = & \sum_{h' \in H} M(h'v, v) = \sum_{u \in V} M(u,v) ,
\end{eqnarray*}
where the second equality follows from automorphism invariance.\vspace*{-2pt}
\end{pf*}
\begin{pf*}{Proof of Proposition \ref{prop:no_minimal_infinite_clusters}}
Consider the mass transport in which each vertex $v$ sitting in a slim
infinite cluster sends unit mass to the unique endpoint of this cluster.
Vertices not sitting in a slim infinite cluster send no mass at all.
Then the expected mass sent from a vertex is at most $1$, while if slim
infinite clusters exist with positive probability then some vertices
will receive infinite mass with positive probability,
so that the expected mass received is infinite, contradicting
(\ref{eq:mass_transport_equilbrium}).
\end{pf*}

Proposition \ref{prop:no_minimal_infinite_clusters} is just an illustrative
example, but BLPS \cite{BLPSgip} proved several other more interesting
reuslts using the Mass-Transport Principle, which turns out
to be quite a potent tool in nonamenable settings where classical density
arguments and ergodic averages
are not available in the same way as in the amenable case.
The Mass-Transport Principle in itself is not especially deep or difficult.
Rather, in the words of BLPS \cite{BLPScritperc}, ``the creative
element in applying the mass-transport method is to make a judicious
choice of the transport function $m(u,v, \omega)$.'' We will see some
examples in this section and the next. (For a remarkable recent
development of the mass-transport method, see Aldous and Lyons \cite{AL}
and Schramm \cite{S08}.)

The following result
characterizes amenability of Cayley graphs (and more generally of
unimodular transitive graphs) in terms of a certain percolation
threshold for invariant percolation.
For a transitive graph $G=(V,E)$ and a~automorphism invariant site
percolation on $G$ with distribution $\mu$ on $\{0,1\}^V$, write
$\pi(\mu)$ for the marginal probability that a given vertex is open.
\begin{theorem}[\cite{BLPSgip}] \label{thm:threshold}
Let $G$ be a unimodular transitive graph, and defi\-ne~$p_{c,\mathrm{inv}}(G)$ is
the infimum over all $p \in[0,1]$ such that any automorphism invariant
site percolation $\mu$ on $G$ with $\pi(\mu) =p$ is guaranteed to produce
at least one infinite cluster with positive probability. Then $p_{c,\mathrm{inv}}(G)<1$
if and only if $G$ is nonamenable.
\end{theorem}

Quantitative estimates for $p_{c,\mathrm{inv}}(G)$ in the nonamenable case are also
provided in \cite{BLPSgip}. Theorem \ref{theorem:threshold_quantitative}
below gives such a bound in the bond percolation case. For
site percolation, BLPS \cite{BLPSgip} show
that if $\pi(\mu) \geq\frac{d(G)}{d(G)+ h(G)}$, where~$d(G)$
is the degree of a vertex in $G$,
and $h(G)$ as before is the isoperimetric constant,
then there is at least one infinite cluster with positive probability. Similar
bonds are given for the quasi-transitive case as well.
For the case $G= \T_d$ the bounds go back to H\"aggstr\"om~\cite{H97},
where they were established using a precursor of the mass-transport method,
and also shown to be sharp. The following bound in the bond percolation case
is in terms of the \textit{edge-isoperimetric constant} $h_E(G)$, defined
by
\[
h_E(G) = \inf_S \frac{|\partial_E S|}{|S|},
\]
where as in Definition \ref{defn:Cheeger} the infimum ranges
over all finite $S \in V$, while
$\partial_ES = \{ \langle u,v \rangle\in E\dvtx u \in S, v\in V
\setminus
S \}$.
Clearly, $h(G) \leq h_E(G) \leq(d(G)-1) h(G)$ when $G$ is transitive,
so such a $G$ is amenable in the sense of Definition \ref{defn:Cheeger}
if and only if $h_E(G)= 0$.
\begin{theorem}[\cite{BLPSgip}] \label{theorem:threshold_quantitative}
Let $G=(V,E)$ be transitive and unimodular,
and consider an automorphism invariant bond
percolation on
$G$ such that for each edge $e \in E$ we have
%
%
\begin{equation} \label{eq:threshold_quantitative}
\Pb(e \mbox{ is open}) \geq\frac{d(G) - h_E(G)}{d(G)} .
\end{equation}
Then the percolation produces an infinite cluster with positive probability.
\end{theorem}

The proof is worth exhibiting here, but in order to be able to
follow the elegant argument from the expository follow-up paper
BLPS \cite{BLPScritperc}, I will be content with considering
the case where (\ref{eq:threshold_quantitative}) holds with strict
inequality.\looseness=1

\begin{pf*}{Proof of (almost) Theorem \ref{theorem:threshold_quantitative}}
For a finite subgraph $G'=(V', E')$ of~$G$,
define its average internal degree
\[
h^*_E(G') = \frac{2|E'|}{|V'|}
\]
and set
%
%
\begin{equation} \label{eq:h-star}
h^*_E(G) = \sup_{G'} h^*_E(G'),
\end{equation}
where the supremum is over all finite subgraphs $G'$ of $G$. For any given
such~$G'$, we have
\[
2|E'| + |\{\langle u,v \rangle\in E\dvtx u \in V', v \in V\setminus V'\}|
\leq d(G) |V'|
\]
with equality if and only if all $e \in E$ with both
endpoints in $V'$ are also in $E'$. Hence, $h^*_E(G)+ h_E(G) = d(G)$, so
the right-hand side of (\ref{eq:threshold_quantitative}) equals
$\frac{h^*_E(G)}{d(G)}$.
Now consider a $\{0,1\}^E$-valued automorphism invariant bond
percolation $X$
on $G$ such that
%
%
\begin{equation} \label{eq:threshold_quantitative_rewritten}
\Pb(e \mbox{ is open}) > \frac{h^*_E(G)}{d(G)}
\end{equation}
for each $e \in E$, and assume for contradiction that it a.s.
produces no infinite cluster. We may define a mass transport
where each vertex
sitting in a finite open cluster counts the number of open edges incident
to it, sends out exactly this amount of mass, and distributes it
equally among
all the vertices sitting in its connected component in the percolation
process. In other words, we take the transport function
to be
\[
m(u,v, \omega) = \cases{
\dfrac{d_\omega(u)}{|K(u)|}, &\quad if $u$ is in a finite component of
$\omega$ and $v \in K(u)$,\vspace*{2pt}\cr
0, &\quad otherwise,}
\]
where $d_\omega(u)$ is the degree of $u$ in $X$, and $K(u)$ is the
set of vertices having an open path in $\omega$ to $u$. Then
(\ref{eq:threshold_quantitative_rewritten}) and the assumption that
$X$ produces no infinite clusters
give
\begin{eqnarray*}
\E\biggl[ \sum_{v \in V} m(u,v,X) \biggr]
& \geq& d(G) \min_{e \in E}\Pb(e \mbox{ is open}) \\[2pt]
& > & d(G) \frac{h^*_E(G)}{d(G)}
= h^*_E(G),
\end{eqnarray*}
while the amount $\sum_{v\in V}m(v,u,X)$ received at $u$ is the
average internal degree of its connected component, which is bounded by
$h^*_E(G)$. Hence,
\[
\E\biggl[ \sum_{v \in V} m(u,v,X) \biggr] >
\E\biggl[ \sum_{v \in V} m(v,u,X) \biggr]
\]
contradicting the Mass-Transport Principle.
\end{pf*}

The next result from BLPS \cite{BLPSgip}
concerns the expected degree of a vertex given that it
belongs to an infinite cluster. Here it seems most natural to consider
the bond percolation case. For $G$ transitive and an invariant
bond percolation on $G$ with distribution $\mu$ that produces
at least one infinite cluster with positive probability,
define $\beta(G, \mu)$
as the expected degree of a vertex given that it belongs
to an infinite cluster, and define $\beta(G)= \inf_\mu\beta(G, \mu)$
where the infimum ranges
over all such automorphism invariant bond percolation processes on
$G$.\vspace*{-3pt}
\begin{theorem}[\cite{BLPSgip}] \label{thm:expected_degree}
For $G=(V,E)$ transitive, we have $\beta(G)=2$ if $G$ is unimodular,
and $\beta(G)<2$ otherwise.\vspace*{-3pt}
\end{theorem}

An exact expression for $\beta(G)$ in the unimodular case is also given,
namely $1 + \inf_{\langle u,v\rangle\in E}\frac{|{\Stab(u)v}|}{|{\Stab(v)u}|}$,
where transitivity implies that the infimum is in fact a minimum. Note how
Theorem \ref{thm:expected_degree} immediately
implies Proposition \ref{prop:no_minimal_infinite_clusters}
above: a slim infinite cluster would have vertices both of degree $1$
and of
degree $2$, and these are the only degrees appearing, so the average
degree would have had to be strictly between $1$ and $2$. That
average degree $2$ is needed is quite intuitive, but what is more surprising
is that it is possible to go below~$2$ in the nonunimodular case. In fact,
the bond percolation on Trofimov's graph in Example
\ref{example:Trofimov_percolation} has $\beta(G, \mu)= \frac{3}{2}$,
and this can be pushed down to $\beta(G, \mu)= \frac{5}{4}$
(which is sharp for Trofimov's graph) by letting
the slim infinite clusters live on grandparent--grandchild rather
than parent--child edges.\looseness=-1

Another striking result in BLPS \cite{BLPSgip} concerns the number of
ends of infinite components: if $G$ is quasi-transitive and unimodular and
$\mu$ is an automorphism invariant site or bond percolation, then
the number of ends of any infinite component must be either $1$, $2$ or
$\infty$. (In the amenable case, the argument of Burton and Keane
\cite{BK}
excludes also the case of infinitely many ends.) Furthermore, for infinite
clusters with infinitely many ends, BLPS \cite{BLPSgip} showed that
that such clusters have expected degree strictly greater than~$2$, and
that they have critical values for site or bond percolation that are
strictly greater than $1$.\looseness=-1

\vspace*{-4pt}\section{No infinite cluster at criticality} \label{sect:at_criticality}

I have yet to mention what is possibly the most striking
result of all from BLPS \cite{BLPSgip}. Namely, this study of
automorphism invariant percolation turned out to have the
following implication for i.i.d. percolation, which is a remarkable step in
the direction of Conjecture \ref{conj:no_perc_at_criticality}.\vadjust{\eject}
\begin{theorem} \label{thm:at_criticality}
Let $G$ be a nonamenable unimodular quasi-transitive graph, and consider
i.i.d. site percolation on $G$ at the critical value $p=p_{c,\mathrm{site}}$.
Then there is a.s. no infinite cluster. The analogous statement for
i.i.d. bond percolation holds as well.
\end{theorem}

Due to the focus in \cite{BLPSgip} being on the more general setting
of automorphism invariant percolation, the proof given
there is not the most direct possible. The authors therefore chose to
publish a separate expository note, BLPS~\cite{BLPScritperc}, with
a more direct proof which is well worth recalling here. Following~\cite{BLPScritperc}, I will restrict to the
case of bond percolation on a (nonamenable, unimodular) transitive
graph; the cases of site percolation and quasi-transitive graphs require
only minor modification.

The reason why unimodularity is needed in the proof is
that the Mass-Transport Principle is used---in fact, it is used
several times
(including in the proof of Theorem \ref{theorem:threshold_quantitative}
which the proof of Theorem \ref{thm:at_criticality} falls back on). But
we should probably expect the result to be true also in the
nonunimodular case (certainly if we trust Conjecture
\ref{conj:no_perc_at_criticality}). See Tim\'ar \cite{Tim}
and Peres, Pete and Scolnicov \cite{PPS} for what are
perhaps the best efforts to date toward a~better understanding of the nonunimodular case.
\begin{pf*}{Proof of Theorem \ref{thm:at_criticality} for bond
percolation on transitive graphs}
Let $G=(V,E)$ be a nonamenable unimodular transitive graph, and
consider an i.i.d. bond percolation $X$ on $G$ with $p=p_{c,\mathrm{bond}}(G)$.
By the Newman--Schulman $0$--$1$--$\infty$ law, the number of infinite
clusters in $X$ is an a.s. constant $N$ which equals either $0$, $1$
or $\infty$, and we need to rule out the possibilities $N=1$ and
$N=\infty$.

\textit{Case} I. \textit{Ruling out $N=1$.}
Assume for contradiction that the percolation
configuration $X \in\{0,1\}^E$ has a unique infinite cluster a.s.
Let $\{Y(e)\}_{e \in E}$ be an i.i.d. collection of random variables,
uniformly distributed on the unit interval $[0,1]$ and independent also
of $X$. For each $\eps\in(0,1)$ and $e \in E$, define\looseness=1
\[
X_\eps(e) = \cases{
1, &\quad if $X(e)=1$ and $Y(e)> \eps$,\cr
0, &\quad otherwise,}
\]\looseness=0
and note that
%
%
\begin{eqnarray} \label{eq:percolation_with_lower_parameter}
X_\eps\in\{0,1\}^E \mbox{ is an i.i.d. bond percolation on $G$}\nonumber\\[-8pt]\\[-8pt]
\eqntext{\mbox{with
parameter } (1- \eps) p_{c, \mathrm{bond}}.}
\end{eqnarray}
For $\eps\in(0,1)$, define yet another bond percolation
$Z_\eps\in\{0,1\}^E$, not i.i.d. but automorphism invariant, as follows.
As before let $\dist_G$ denote graph-theoretic distance
in $G$, and for each $v \in V$ define $U(v)$ as the set of vertices in
the infinite\vadjust{\eject} cluster of $X$ that minimize $\dist_G(u,v)$. Note
that $U(v)$ is finite for all $v \in V$. For
an edge $e = \langle v,w \rangle\in E$, set
\[
Z_\eps(e) =
\cases{1, &\quad if all vertices in $U(v)$ and $U(w)$ are\cr
&\quad in the same connected component of $X_\eps$,\cr
0, &\quad otherwise.}
\]
This defines the percolation $Z_\eps\in\{0,1\}^E$. For any
$\langle v,w \rangle\in E$, there exists some
finite collection $T(v,w) \subset E$
of open edges in $X$ that together connect all the vertices in
$U(v)$ and $U(w)$ to each other (this is where the assumption $N=1$ is used).
For definiteness, we take $T(v,w)$ to be the edge set with minimal
cardinality having this property, and with minimization of
$\sum_{e\in T(v,w)}Y(e)$ acting as tie-breaker.
Each edge $e'$ in this collection has $Y(e')>0$ a.s., and so
$\min_{e' \in T(v,w)} Y(e') >0$ so that $\lim_{\eps\rightarrow
0}Z_\eps
(e)=1$ a.s.
Hence,
%
%
\begin{equation} \label{eq:sending_eps_to_0}
\lim_{\eps\rightarrow0} \Pb[Z_\eps(e)=1 ] = 1 ,
\end{equation}
and Theorem \ref{theorem:threshold_quantitative} ensures that that
the percolation process $Z_\eps$ contains an infinite cluster
with positive probability. But when $Z_\eps$ contains an
infinite cluster, then so does $X_\eps$. In view of
(\ref{eq:percolation_with_lower_parameter}), this contradicts the definition
of~$p_{c,\mathrm{bond}}$, so we are done with Case I.

\textit{Case} II. \textit{Ruling out $N=\infty$.} This time assume for
contradiction that~$X$ contains infinitely many infinite clusters a.s. Here we need the
concept of \textit{encounter points} introduced by Burton and Keane
\cite{BK}
in their famous short proof of uniqueness of the infinite cluster on
$\Z^d$.
An encounter point in a~percolation process $X$ is a vertex $v \in V$ that has three disjoint
open paths to infinity that would fall in different connected
components of $X$ if the vertex $v$ were to be removed, but does
not have four such paths. Burton
and Keane showed that if $N=\infty$, then $X$ contains encounter
points a.s.; their proof was formulated for $G=\Z^d$, but goes
through unchanged for general (quasi-)transitive graphs.

BLPS \cite{BLPScritperc} begin by noting that
%
%
\begin{equation} \label{eq:encounter_points_everywhere}
\begin{tabular}{p{325pt}}
if $v \in V$ is an encounter point, then a.s. each of the three
infinite clusters $C_1(V), C_2(v), C_3(v)$ that the removal
of $v$ would produce contains further encounter points.
\end{tabular}\hspace*{-20pt}
\end{equation}
To see this, consider the mass transport in which each vertex $u$ sitting
in an infinite cluster with encounter points sends
unit mass to the nearest encounter point with respect to $\dist_X$,
splitting it equally in case of a tie; here $\dist_X$ means
graph-theoretic distance in the open subgraph of $G$ defined by $X$. Failure
of (\ref{eq:encounter_points_everywhere}) would cause a contradiction
to the Mass-Transport Principle similarly as in the proof of Proposition
\ref{prop:no_minimal_infinite_clusters}.

Next, we go on to define a random graph $H=(W,F)$, whose vertex set
$W \subset V$ is the set of encounter points in $X$, and whose edge set
$F$, which we are about to specify, will not necessarily be a subset of $E$
(so $H$ is not a subgraph of $G$). Let $\{Y(v)\}_{v \in W}$ be i.i.d.,
uniformly distributed on $[0,1]$ and independent of $X$. Each $v \in W$
selects three other $u_1, u_2, u_3 \in V$ to form edges to, according
to the following rule: one $u_i$ should be chosen in each of the components
$C_1(v)$, $C_2(v)$ and $C_3(v)$, and in each such component~$u_i$ is chosen
to minimize $\dist_X(v, u_i)$, with minimization of $Y(u_i)$ acting as
a tie-breaker. An equivalent way to formulate this is that $u_i$ is chosen
in $C_i$ to minimize the ``distance'' $\dist_{X,Y}(v, u_i)$ defined
as $\dist_{X,Y}(v, u_i)= \dist_{X}(v, u_i)+ Y(v) + Y(u_i)$.

Each $v \in V$ thus gets $H$-degree at least $3$, but the $H$-degree may
exceed~$3$ if $v$ is selected by some $w \in W$ which is not among $v$'s
preferred triplet. An application of the Mass-Transport Principle shows that
the expected number of vertices that choose $v$ is exactly $3$,
so the expected $H$-degree of $v$ (conditional on being in $W$)
is somewhere between $3$ and $6$, and in particular its
degree is a.s. finite.

A crucial step of the argument is now to show that
%
%
\begin{equation} \label{eq:no_cycles}
\mbox{the graph $H$ has no cycles.}
\end{equation}
To see this, we first note that if $v \in W$ is in such a cycle,
then its two neighbors in this cycle must belong to the same $C_i(v)$,
as otherwise we would get a direct contradiction to the definition
of an encounter point. Using this, it is not hard to see (or consult
\cite{BLPScritperc} for a more detailed argument)
that any cycle
$v_1 \leftrightarrow v_2 \leftrightarrow
\cdots\leftrightarrow v_k \leftrightarrow v_1$ would
have to satisfy either
\[
\dist_{X,Y}(v_1, v_2) < \dist_{X,Y}(v_2, v_3)
< \cdots< \dist_{X,Y}(v_k, v_1) < \dist_{X,Y}(v_1, v_2)
\]
or
\[
\dist_{X,Y}(v_1, v_2) > \dist_{X,Y}(v_2, v_3)
> \cdots> \dist_{X,Y}(v_k, v_1) > \dist_{X,Y}(v_1, v_2),
\]
which in either case is of course a contradiction. Hence,
(\ref{eq:no_cycles}).

Now take an $\eps>0$ and consider as in Case I the $\eps$-thinned
percolation process $X_\eps\in\{0,1\}^E$. Using $X_\eps$, we define
the subgraph $H_\eps=(W, F_\eps)$ obtained from $H$ by deleting
each $e = \langle v,w\rangle\in F$
such that $v$ and $w$ fail to be in the same connected component of~$X_\eps$. By the definition of
$p_{c,\mathrm{bond}}$, we have that~$X_\eps$ has
no infinite clusters, so that~$H_\eps$ has no infinite clusters
either.

For $v \in W$, write $K_\eps(v)$ for the set of vertices in $W$ that belong
to the same connected component of $H_\eps$ as $v$. Also write
$\partial_{\mathrm{int}}K_\eps(v)$ ($\mathrm{int}$ as in ``internal boundary'') for the
set of vertices in $K_\eps(v)$ that have at least one neighbor in $H$ which
is not in $K_\eps(v)$.

Now define the following mass transport on $G=(V,E)$. No vertex \mbox{$v\in V$} sends
any mass unless it is an encounter point, that is, unless $v \in W$. Each
encounter point $v$ sends unit mass and divides it equally amongst the
vertices in $K_\eps(v)$. This defines the mass transport, and
the mass received at $v$ becomes
\[
\cases{
\dfrac{|K_\eps(v)|}{|\partial_{\mathrm{int}}K_\eps(v)|}, &\quad
if $v$ is an encounter point and belongs to
$\partial_{\mathrm{int}}K_\eps(v)$,\vspace*{2pt}\cr
0, &\quad otherwise.}
\]
Since $F$ is a forest in which each vertex has
degree is at least $3$, its isoperimetric
constant is easily seen to be at least $1$ (which holds with
equality on the binary tree~$\T_2$), and similarly
$|K_\eps(v)|/|\partial_{\mathrm{int}}K_\eps(v)| \leq2$. So the
expected mass received at $v$ is bounded by
$2\Pb(v \in W, v \in\partial_{\mathrm{int}}K_\eps(v))$,
while the expected mass sent from $v \in V$ equals $\Pb(v \in W)$.
Hence, the Mass-Transport Principle gives
%
%
\begin{equation} \label{eq:yet_another_mass_balance}
\Pb(v \in W) \leq2 \Pb\bigl(v \in W, v \in\partial_{\mathrm{int}}K_\eps(v)\bigr).
\end{equation}
But similarly as in the argument for (\ref{eq:sending_eps_to_0}), we get
for any $v \in W$ that $v \notin\partial_{\mathrm{int}}K_\eps(v)$ for all
sufficiently small $\eps$, so that
\[
\lim_{\eps\rightarrow0} 2 \Pb\bigl(v \in W, v \in\partial_{\mathrm{int}}K_\eps(v)\bigr)
= 0 .
\]
Since (\ref{eq:yet_another_mass_balance}) was shown to hold
for any $\eps>0$, we get $\Pb(v \in W)=0$, which contradicts
$N=\infty$,
so Case II is finished and the proof is complete.
\end{pf*}

\section{Random walks on percolation clusters}
\label{sect:random_walk}

One of the most natural probabilistic objects, besides i.i.d. percolation,
to define on a graph $G=(V,E)$, is simple random walk (SRW), which is a
$V$-valued random process $\{Z_0, Z_1, \ldots\}$ where
one (typically) takes $Z_0=\rho$ for some prespecified choice of
$\rho\in V$, and then iterates the following: given
$Z_0, Z_1, \ldots, Z_{n-1}$, the value of $Z_n$ is chosen uniformly
among the neighbors in $G$ of $Z_{n-1}$.

In contrast to the study of percolation on nonamenable Cayley graphs
and related classes of graphs which
began to take off only in the 1990s, the literature
on SRWs on such graphs goes back much further. An early seminal
contribution is the work of Kesten \cite{K59,K59b}
from the late 1950s showing that if
$G$ is a Cayley graph for a finitely generated group, then the return
probability $\Pb[Z_n=\rho]$ decays exponentially if and only if $G$
is nonamenable. See, for example, Woess \cite{Woess} for an
introduction to
this field.

Another topic of considerable interest is the study of SRW
on percolation clusters. For the $\Z^d$ case, see, for instance, papers
like \cite{DMFGW} and \cite{BerBisk} on central limit theorems,
and \cite{GKZ} which extends P\'olya's classical $d=2$ versus
\mbox{$d\geq3$} recurrence-transience dichotomy for random walk on $\Z^d$
to the case of supercritical percolation on $\Z^d$.

Given these traditions, it was a very natural step for Schramm and
his collaborators to go on to consider random walks on
percolation clusters on nonamenable graphs. Their work is of two kinds:
on one hand, the analysis of SRW on a percolation cluster as a
worthwhile object of study in its own right, and, on the other hand,
the exploitation of random walk on a percolation cluster as a means
toward understanding properties of percolation clusters that do not
primarily have anything to do with random walk. A remarkable application
of the second kind will be described in
Section \ref{sect:cluster_indistinguishability} on so-called
cluster indistinguishability, while in the present section I will
recall a result of the first kind (Theorem \ref{thm:positive_speed} below)
from a rich paper by Benjamini, Lyons
and Schramm, henceforth BLS \cite{BLS}.

A natural question to ask for SRW on an infinite graph $G$ is how
fast it escapes from the starting point $\rho$, that is, how fast
does $\dist_G(\rho, Z_n)$ grow? Define the \textit{speed}
\[
S = \lim_{n \rightarrow\infty} \frac{\dist_G(\rho, Z_n)}{n}
\]
provided the limit exists. When $G$ is the $\Z^d$ lattice,
$\dist_G(\rho, Z_n)$ scales like~$\sqrt{n}$, and not surprisingly
$S=0$ a.s. More generally, when $G$ is any transitive graph,
the Subadditive Ergodic Theorem immediately
implies that the limit $S$ exists and
is an a.s. constant.

If we go on on to consider the speed of SRW on an infinite cluster
of, say, i.i.d. bond percolation on $G$ (still with
the speed defined with respect to $\dist_G$),
then the existence of the speed $S$
is less obvious. However, BLS \cite{BLS} showed, when $G$ is unimodular,
and the percolation process is automorphism invariant,
that the speed does exist a.s. and does not depend on the random walk, but
only on the percolation configuration. Having come that far, it is easy to
see that the speed cannot depend on where in an infinite cluster the
random walk starts, so each infinite cluster has a well-defined characteristic
SRW speed. For the i.i.d. bond percolation case, we can then invoke
the cluster indistinguishability result of Lyons and Schramm
\cite{LS} (Theorem \ref{thm:cluster_indistinguishability} below) to
deduce that all infinite clusters have the same SRW speed.

Of particular interest is to determine whether the speed is
zero or positive. For the nonamenable unimodular case,
BLS \cite{BLS} found a general answer.
\begin{theorem} \label{thm:positive_speed}
The speed $S$ of SRW on an infinite cluster of an i.i.d. bond percolation
$X$ on a unimodular nonamenable transitive graph $G$ satisfies $S>0$ a.s.
\end{theorem}

An important step in the proof of this is the following result (also of
independent interest) from \cite{BLS} on the geometry of infinite
clusters. That the nonamenability of $G$ should be inherited by
the infinite clusters of $X$ is too much to hope for: a sufficient condition
for a percolation cluster in $X$ to be amenable is that it contains
arbitrarily long ``naked'' paths, that is, paths of vertices with
$X$-degree $2$, and it can be shown that a.s. all infinite clusters arising
from i.i.d. percolation on a transitive graph contain such paths. But the
infinite clusters of $X$ do contain nonamenable subgraphs.
\begin{theorem} \label{thm:nonamenable_subgraph}
Any infinite cluster in i.i.d. bond percolation
$X$ on a~uni\-modular nonamenable transitive graph $G$ contains
a nonamenable subgraph a.s.\vspace*{-2pt}
\end{theorem}

The proof in \cite{BLS} of this result involves yet another
application of mass transport. Both Theorems \ref{thm:nonamenable_subgraph}
and \ref{thm:positive_speed} are in fact proved more generally
than just for i.i.d. percolation. Automorphism invariance alone does
not suffice (each of Examples 3.1, 3.2 and 3.3 of
H\"aggstr\"om \cite{H97} shows this, and moreover that the $>$ in
condition (d) below cannot be replaced by a $\geq$), but if we add
any of the conditions:
\begin{enumerate}[(a)]
\item[(a)] $X$ is i.i.d.,
\item[(b)] $X$ has a unique infinite cluster a.s.,
\item[(c)] the infinite clusters of $X$ have at least $3$ (hence
infinitely many) ends a.s., or
\item[(d)] $X$ is ergodic with sufficiently large values of
$\Pb(e \mbox{ is open})$, more precisely the expected degree of
a vertex should strictly exceed the
quantity $h^*_E$ defined in~(\ref{eq:h-star}),
\end{enumerate}
then the conclusions of Theorems \ref{thm:nonamenable_subgraph}
and \ref{thm:positive_speed} hold; cf. Theorems 3.9 and 4.4 in BLS
\cite{BLS}.

\vspace*{-2pt}\section{Cluster indistinguishability}
\label{sect:cluster_indistinguishability}

Consider an i.i.d. bond percolation $X$ with parameter $p>p_{c,\mathrm{bond}}(C)$
on a graph $G$, so that $X$ produces
one or more infinite clusters. It is then natural to ask questions about
properties of these infinite clusters. Properties that
we have already discussed in earlier sections
include the number of ends of an infinite cluster,
whether it contains encounter points, and the speed of SRW on the cluster.

Yet another natural such property of an infinite cluster $C$ is
the value of $p_{c,\mathrm{bond}}(C)$, that is, how much further i.i.d.
edge-thinning can the infinite cluster $C$ take before it breaks apart
into finite components only? For $G$ transitive, it is known (see
H\"aggstr\"om and Peres \cite{HP} for the unimodular case,
and Schonmann \cite{Scho} for the full result) that a.s.
$p_{c,\mathrm{bond}}(C)= p_{c,\mathrm{bond}}(G)/p$ for all infinite clusters $C$ of $X$.

All these properties are examples of \textit{invariant properties}. For
$G=(V,E)$ transitive with automorphism group $\Aut(G)$, a property (which
may or may not hold for clusters of bond percolation on $G$)
can be identified with a Borel measurable subset
of $\{0,1\}^E$, and a property $\cA\subset\{0,1\}^E$ is said to be
invariant if for all $\omega\in\cA$ and all $\gamma\in\Aut(G)$
we have $\gamma\omega\in\cA$. Lyons and Schramm~\cite{LS} proved
the following Theorem \ref{thm:cluster_indistinguishability},
known as \textit{cluster indistinguishability}. Shortly before \cite{LS},
a~weaker result for
the case of so-called increasing invariant properties was established
in \cite{HP}.
\begin{theorem} \label{thm:cluster_indistinguishability}
Let $G=(V,E)$ be a nonamenable unimodular transitive graph, and
consider i.i.d. bond percolation
$X$ on $G$ with $p$ in the parameter regime where $X$ produces infinitely
many infinite clusters a.s.
Then, for any invariant component property $\cA$, we have
a.s. that either all infinite components of $X$ satisfy $\cA$ or
all infinite components of $X$ satisfy $\neg\cA$.
\end{theorem}

Space does not permit me to give the full proof of this beautiful result,
but I can explain what the main steps are.
\begin{pf*}{Sketch proof of Theorem \ref{thm:cluster_indistinguishability}}
Let $G=(V,E)$ and the percolation $X\in\{0,1\}^E$ be as in the theorem,
and assume for contradiction that $\cA$ is an invariant property such
that with positive probability, $X$ contains both infinite clusters
with property $\cA$ and infinite clusters with property $\neg\cA$.

\textit{Step} I. \textit{Existence of pivotal edges.}
For an infinite cluster $C$ of $X$ and an edge $e\in E$ with $X(e)=0$
that has an endpoint in $C$, call $e$ \textit{pivotal} for $C$ if either
$C \in\cA$ and switching on the edge $e$ would create an infinite cluster
(containing $C$) with property $\neg\cA$, or vice versa. If there exists
an $e\in E$ with $X(e)=0$ that has one endpoint in an infinite cluster
with property $\cA$ and the other in an infinite cluster
with property $\neg\cA$, then clearly $e$ is pivotal for one of the clusters.
Otherwise, there exists such a pair of clusters within finite distance
from each other, and by sequentially switching on one edge after another
on a finite path between them we see that somewhere along the way
one infinite cluster of type $\cA$ must turn into $\neg\cA$ or vice versa.
(This is an example of a well-known technique in percolation theory known
as \textit{local modification}, pioneered by Newman and Schulman \cite
{NS} and
Burton and Keane \cite{BK}; cf. also Coupling 2.5 in \cite{HJ06} for
a careful explanation.) Hence, pivotal edges exist with positive probability.

\textit{Step} II. \textit{Stationarity of random walk.}
Consider a SRW $\{Z_0, Z_1, \ldots\}$
on a~percolation cluster of $X$, defined as in Section
\ref{sect:random_walk}. It would be nice to think that
the percolation configuration
``as seen from the point of view of the walker,''
would be stationary. To make this precise, fix for each $v \in V$
a~$\gamma_{v, \rho} \in\Aut(G)$ that maps $v$ on $\rho$, where, as
in Section \ref{sect:random_walk}, $\rho$ is the starting point of the
random walk. The idea of stationarity of the percolation configuration
as seen from the random walker is that
for any Borel measurable $\Stab(\rho)$-invariant
$\cB\in\{0,1\}^E$ and any $n$
we should have
\[
\Pb(\gamma_{Z_n, \rho}X \in\cB) = \Pb(X \in\cB) .
\]
This is not
true, however, and to see this, think, for example, about what happens
when the
starting point $\rho$ happens to be in a finite
open cluster~$C$ which is simply
a path of length $2$; in other words, $C$ has three vertices, one of
which has degree $2$ and two of which have degree $1$. Conditioned
on $\rho$ being in such a cluster, it has probability $\frac{1}{3}$
of being in the vertex of degree~$2$ (this follows from a straightforward
mass-transport argument). But SRW on such an open cluster
will spend half the time (either all even or all odd times) on that
vertex, so it cannot possibly be stationary in the desired
sense.\looseness=-1

All is not lost, however. Stationarity can be recovered by a minor
modification of SRW, namely the \textit{delayed simple random walk}
(\textit{DSRW}),
denoted $\{\tilde{Z}_0, \tilde{Z}_1, \ldots\}$.
This is again a $V$-valued random process, and as with SRW
we take $\tilde{Z}_0= \rho$, but the transition mechanism is slightly
different: given $\tilde{Z}_0, \ldots, \tilde{Z}_{n-1}$, a vertex $w$
is chosen according to uniform distribution on the set of $\tilde{Z}_{n-1}$'s
$G$-neighbors, and we set
\[
\tilde{Z}_n = \cases{
w, &\quad if $X(\langle\tilde{Z}_{n-1}, w \rangle)=1$, \cr
\tilde{Z}_{n-1}, &\quad otherwise.}
\]
Such a DSRW has the desired stationarity property, that is,
\[
\Pb(\gamma_{\tilde{Z}_n, \rho}X \in\cB) = \Pb(X \in\cB)
\]
for any $\cB$ and any $n$, and this can be shown via yet
another mass-transport argument.

(Stationarity of DSRW is clearly an interesting result in its own right,
and the idea was first exploited in \cite{H97}.
Lyons and Schramm decided to put their most general
version of the stationarity result
not in \cite{LS} but in a separate paper \cite{LS99b}.
Without unimodularity, however, stationarity fails, as
is easily seen by considering DSRW in
Example \ref{example:Trofimov_percolation}: here, $\tilde{Z}_0$ has
probability $\frac{1}{2}$ of sitting in the ``topmost'' (closest to
$\xi$)
component of its open cluster, while the probability
that $\tilde{Z}_n$ does so tends to $0$ as $n \rightarrow\infty$.)

\textit{Step} III. \textit{Lots of encounter points.} Since $X$ has infinitely many
infinite clusters, we can modify $X$ by changing finitely
many edges so as to connect three of them to form an encounter point.
Hence, encounter points exist with positive probability; this is just the
local modification argument of Burton and Keane \cite{BK}.
By (\ref{eq:encounter_points_everywhere}), we have that every infinite
cluster with encounter points must in fact have infinitely many.
From this, it follows that \textit{every} infinite cluster has infinitely
many encounter points, because otherwise we could use the local modification
technique to connect an infinite cluster with encounter points
to one without and obtain a contradiction to
(\ref{eq:encounter_points_everywhere}).

\textit{Step} IV. \textit{Random walk is transient.}
Given the prevalence of encounter points, infinite clusters contain,
loosely speaking, a large-scale
structure similar to the binary tree $\T_2$. This strongly suggests
that DSRW on such an infinite cluster should be transient.
Lyons and Schramm \cite{LS} converts this intuition into a proof by
combining the result mentioned in the final paragraph of
Section \ref{sect:BLPS-gip} about such infinite clusters $C$
having $p_{c,\mathrm{bond}}(C)<1$,
with a basic comparison of random walk and percolation on trees due to
Lyons~\cite{L90}. (Alternatively, we could quote
Theorem~\ref{thm:positive_speed} here, but that would be a~detour.)

\textit{Step} V. \textit{Local modification applied to a pivotal edge.}
Let $\cB$ be the event that the starting point
$\rho$ of the DSRW is in an infinite cluster with property~$\cA$.
Given some small
$\eps>0$, we can find a $k< \infty$ and a $\Stab(\rho)$-invariant event
$\cB^*$ that depends
only on edges within $G$-distance $k$ from $\rho$, approximating
$\cB$ in the sense that
%
%
\begin{equation} \label{eq:approximate_property}
\Pb(\cB\triangle\cB^*)< \eps.
\end{equation}
Define, for each $n$,
%
%
\begin{equation} \label{eq:def_of_Y}
Y_n = \cases{
1, &\quad if $\gamma_{\tilde{Z}_n,\rho}X \in\cB^*$, \cr
0, &\quad otherwise.}
\end{equation}
Then $\{Y_0, Y_1, \ldots\}$ is a stationary process, so the limit
$\bar{Y}=\lim_{n\rightarrow\infty}\frac{1}{n}\sum_{i=0}^{n-1}Y_n$
exists a.s.
Furthermore,
%
%
\begin{equation} \label{eq:depends_only_on_percolation}
\begin{tabular}{p{325pt}}
the limit $\bar{Y}$ depends only on the percolation
configuration and not on the DSRW.
\end{tabular}\hspace*{-25pt}
\end{equation}
To see (\ref{eq:depends_only_on_percolation}), define a second
DSRW $\{\tilde{Z}'_0, \tilde{Z}'_1, \ldots\}$ from $\rho$ which
conditionally
on $X$ is independent of $\{\tilde{Z}_0, \tilde{Z}_1, \ldots\}$,
define $\{Y'_0, Y'_1, \ldots\}$ analogously as in~(\ref{eq:def_of_Y}),
and set $\bar{Y}'=\lim_{n\rightarrow\infty}\frac{1}{n}\sum
_{i=0}^{n-1}Y'_n$.
Then we can define $Y_{-i}=Y'_i$ for each $i\geq1$, and it turns
out (see \cite{LS}, Lemma 3.13) that the two-sided sequence
$\{\ldots, Y_{-1}, Y_0, Y_1, \ldots\}$ becomes stationary. Now, if
(\ref{eq:depends_only_on_percolation}) failed we would with positive
probability have $\bar{Y}\neq\bar{Y}'$, which however would contradict
stationarity of the two-sided sequence $\{\ldots, Y_{-1}, Y_0, Y_1,
\ldots\}$
in view of the ergodic theorem. Hence, (\ref{eq:depends_only_on_percolation}).

Thus, we can assign each infinite cluster $C$ a value $\alpha(C)$ as
the a.s. value of $\bar{Y}$ that DSRW on that cluster would produce.
By (\ref{eq:approximate_property}), we have that infinite clusters
with different values of $\alpha(C)$ coexist with positive probability,
provided we chose $\eps$ small enough to begin with. Call
an edge $e\in E$ $\alpha$-\textit{pivotal}
if it is closed [$X(e)=0$] with endpoints
in two different infinite clusters $C$ and $C'$ with
%
%
\begin{equation} \label{eq:alpha-pivotality}
\alpha(C) \neq\alpha(C') .
\end{equation}
By the reasoning in Step 1, there exist $\alpha$-pivotal edges
in $X$
with positive probability. In particular, there exists an $e \in E$
which with positive probability is $\alpha$-pivotal for the infinite
cluster containing $\rho$. Suppose this happens, let $C$ be the infinite
cluster containing $\rho$, and let $C'$ be the other infinite cluster
meeting $e$. Then $\bar{Y}=\alpha(C)$ a.s.

But look what happens if we apply local modification
to the edge $e$. If we turn $e$ on
[$X(e)=1$] and leave the status of all other edges intact, we get
a~new infinite cluster $C''$ uniting $C$ and $C'$. What will then happen
with $\bar{Y}$? By transience of DSRW, we get with positive probability
that the DSRW escapes to infinity without ever noticing the edge $e$,
and we get a.s. on this event that $\bar{Y}=\alpha(C)$; this
uses the fact that $Y_n$ is a function of edges within bounded
distance $k$ from $\tilde{Z}_n$ only. On the other hand,
we get with positive probability that the DSRW reaches $e$, crosses it,
and then escapes to infinity without ever crossing it back; on this event
we get a.s. $\bar{Y}=\alpha(C')$. In view of (\ref{eq:alpha-pivotality}),
this contradicts (\ref{eq:depends_only_on_percolation}) and proves the
theorem.
\end{pf*}

An inspection of the proof to determine which properties of i.i.d.
percolation are actually used reveals that it is enough
to assume automorphism invariance plus
so-called \textit{insertion tolerance}, which is the term
Lyons and Schramm \cite{LS} used for a refinement
of the finite energy property considered by
Newman and Schulman~\cite{NS}, Burton and Keane \cite{BK}, and
others.\looseness=1
\begin{defn} \label{defn:insertion_tolerance}
A bond percolation $X$ on a graph $G=(V,E)$
is called \textit{insertion tolerant} if it admits conditional probabilities
such that for every $e\in E$ and every $\xi\in\{0,1\}^{E \setminus\{
e\}}$
we have $\Pb(X(e)=1 | X(E\setminus\{e\} = \xi)>0$. If instead
$\Pb(X(e)=0 | X(E\setminus\{e\} = \xi)>0$ for all such $e $ and $\xi$,
then $X$ is called \textit{deletion tolerant}.
\end{defn}

(The conjunction of insertion tolerance and deletion tolerance is precisely
the finite energy property.)
\begin{theorem} \label{thm:cluster_indistinguishability_IT}
Let $G=(V,E)$ be a nonamenable unimodular transitive graph, and
consider an insertion tolerant automorphism invariant bond percolation
$X$ on $G$. Then, for any invariant component property $\mathcal{A}$,
we have
a.s. that either all infinite components of $X$ are in $\mathcal{A}$ or
no infinite components of $X$ are in~$\mathcal{A}$.
\end{theorem}

In fact, the formulation in \cite{LS} is even more general than this: here,
as in much of the work discussed in Section \ref{sect:BLPS-gip},
invariance under the full automorphism group can be weakened to invariance
under certain subgroups. Also, the result holds with site percolation
in place of bond percolation.

It is worth mentioning some limitations to the scope of cluster
indistinguishability. For instance, insertion tolerance cannot be
replaced by deletion tolerance in Theorem
\ref{thm:cluster_indistinguishability_IT}, as the following
example from \cite{LS} shows.
\begin{example}
Let $G=(V,E)$ be the binary tree $\T_2$,
recall that\break $p_{c,\mathrm{bond}}(T_2)=\frac{1}{2}$,
and fix $p_1$ and $p_2$ such that $\frac{1}{2}<p_1<p_2<1$. Let
$X \in\{0,1\}^E$ be i.i.d. bond percolation on $G$ with parameter $p_2$.
Then obtain another bond percolation process $X' \in\{0,1\}^E$
from $X$ as follows.
For each infinite cluster~$C$ of $X$ independently, toss a fair coin.
If the
coin comes up heads, delete each edge of $C$ independently with probability
$1 - \frac{p_1}{p_2}$; otherwise let all of $C$'s edges be intact. While
$X'$ is automorphism invariant and deletion tolerant, some of its infinite
clusters $C'$ will have all the characteristics of those produced by
i.i.d. $({p_2})$ percolation and thus have $p_{c,\mathrm{bond}}(C')=\frac{1}{2p_2}$,
while others will look like i.i.d. $({p_1})$ percolation clusters and
thus have
$p_{c,\mathrm{bond}}(C')=\frac{1}{2p_1}$, so cluster indistinguishability fails for
$X'$.
\end{example}

How about the unimodularity assumption in Theorems
\ref{thm:cluster_indistinguishability} and
\ref{thm:cluster_indistinguishability_IT}? This cannot be dropped, not even
from Theorem \ref{thm:cluster_indistinguishability}. To see this, let
$G=(V,E)$ be Trofimov's graph (Example \ref{ex:Trofimov}), and consider
i.i.d. bond percolation with parameter $p \in(p_{c,\mathrm{bond}}(G), 1)$.
Each infinite cluster $C$ then has a ``topmost'' vertex $w(C)$ (in the
direction of the designated end $\xi$), and for $i=0,1,2$ we may define
the property $\cA_i$ by stipulating that $C\in\cA_i$ if $w(C)$ is
directly linked to exactly $i$ of its $\xi$-children via open edges.
Then $\cA_0$, $\cA_1$ and $\cA_2$ are $\Aut(G)$-invariant properties,
and the percolation will a.s. produce infinite clusters of all three kinds,
so cluster indistinguishability fails.

\vspace*{-2pt}\section{In the hyperbolic plane} \label{sect:planar}

In their \textit{Percolation beyond $\Z^d$} paper \cite{BS96}, Benjamini
and Schramm emphasized three properties of graphs that could be expected
to be of particular relevance to the behavior of percolation processes:
quasi-transitivity, (non-)amenability, and planarity. The first two
I have discussed at some length, while the third was only mentioned in
passing in Section \ref{sect:BS-beyond}. In this section, I will briefly
make up for this.

Planarity plays a crucial role in the classical study of percolation
on the~$\Z^2$ lattice, such as in the seminal contributions
by Harris~\cite{Har} and Kesten~\cite{K}. A~key device is the notion of
planar duality: for a graph $G=(V,E)$ with a~planar embedding in $\R^2$
(or in some other two-dimensional manifold), it is often
useful to define its dual graph $G^\dagger=(V^\dagger,E^\dagger)$ by
identifying $V^\dagger$ with the faces of the planar embedding of $G$,
and including an edge $e^\dagger\in E^\dagger$ crossing each $e\in E$.
If $X\in\{0,1\}^E$ a bond percolation on $G$, then we can define a~bond
percolation $X^\dagger\in\{0,1\}^E$ on $G^\dagger$ by declaring,
for each $e \in E$, $X^\dagger(e^\dagger)=1-X(e)$. If $X$ is i.i.d. ($p$),
then $X^\dagger$ becomes i.i.d. ($1-p$). Furthermore,
if $G$ is the $\Z^2$ lattice,
then $G^\dagger$ is isomorphic to $G$, and for $p=1-p=\frac{1}{2}$ the
distributions of $X$ and $X^\dagger$ will be the same; this observation
is basic to proving the Harris--Kesten theorem
that $p_{c,\mathrm{bond}}(\Z^2)=\frac{1}{2}$.

Suppose that $G$ is infinite, planar and transitive.
It is known (see Babai~\cite{Bab}) that such a graph, equipped with
the usual distance $\dist_G$, is quasi-isometric
to exactly one of the four spaces $\R$, $\R^2$, $\T_2$ and the
hyperbolic plane~$\Hb^2$, the last of which may be defined as the open unit
disk $\{z \in\C\dvtx|z|<1\}$ equipped with the metric $s$
given by
\[
ds^2 = \frac{dx^2 + dy^2}{(1-x^2-y^2)^2} .
\]
In the paper \textit{Percolation in the hyperbolic plane} \cite{BS01},
Benjamini and Schramm consider the situation where $G$ has one end, in
which case it must be quasi-isometric to either $\R^2$ or $\Hb^2$.
Which of these is determined by amenability: if $G$ is amenable, then
$\R^2$, while if it is nonamenable, then $\Hb^2$. The focus of \cite{BS01}
is on the nonamenable case; hence, the title of the paper. A happy
circumstance here is that a planar nonamenable transitive graph with one
end is also unimodular (Proposition 2.1 in \cite{BS01}), which
allows the machinery developed in BLPS \cite{BLPSgip} to come into
play.\vadjust{\eject}

Recall the Newman--Schulman $0$--$1$--$\infty$ law about the number
of infinite clusters on transitive graphs. This result needs the percolation
process to be automorphism invariant and insertion tolerant
(cf. Definition \ref{defn:insertion_tolerance}). In general, insertion
tolerance cannot be dropped (not even on the $\Z^2$ lattice;
see Burton and Keane \cite{BK91}), although in the present setting,
remarkably enough, it can:
\begin{theorem}[(Theorem 8.1 in \cite{BLPSgip})]
\label{thm:Newman-Schulman_without_insertion_tolerance}
If $G=(V,E)$ is a planar nonamenable quasi-transitive graph with one end,
then a.s. any
automorphism invariant bond percolation $X \in\{0,1\}^E$ has $0$, $1$,
or infinitely many infinite clusters.
\end{theorem}
\begin{pf}
This proof from \cite{BS01} differs somewhat
from the original one in~\cite{BLPSgip}. There is no loss of generality
in assuming that the number of infinite clusters is an a.s. constant $k$.
Assume for contradiction that $k\in\{2,3,\ldots\}$. By randomly deleting
all edges of all infinite clusters but two, chosen uniformly at
random, we still preserve automorphism invariance; thus we may assume
$k=2$. Call one of them $C_1$ and the other $C_2$, using a fair coin
toss to decide which is which. Then turn on each edge $e$
which does not meet $C_2$ and from which
there is a path in $G$ to $C_1$ that does not meet $C_2$.
This preserves automorphism invariance, and expands $C_1$ to a larger
infinite cluster $\hat{C}_1$ that, loosely speaking, sits as close
to $C_2$ as is possible without touching it.
Now define a bond percolation $\hat{X}^\dagger$ on the
dual graph $G^\dagger$ (which is also planar, nonamenable, quasi-transitive
and one-ended) by turning on exactly those\vspace*{2pt} edges $e^\dagger\in
E^\dagger$
whose corresponding $e \in E$ are pivotal for connecting
$\hat{C}_1$ to $C_2$. Then $\hat{X}^\dagger$ is $\Aut(G^\dagger
)$-invariant,
and consists (due to one-endedness of~$G$) of a single bi-infinite open path.
Hence, $\hat{X}^\dagger$ has a unique infinite cluster~$C^\dagger$ with
$p_{c,\mathrm{bond}}(C^\dagger)=1$. On the other hand, the reasoning in the
proof of Theorem \ref{thm:at_criticality}, Case I, shows that
any unique infinite cluster arising from automorphism invariant percolation
in such a graph has $p_{c,\mathrm{bond}}<1$,
and this is the desired contradiction.
\end{pf}

From here, Benjamini and Schramm \cite{BS01} go on to show a number
of interesting results for percolation\vadjust{\goodbreak} in the hyperbolic plane. For starters,
let $G$ be as
in Theorem \ref{thm:Newman-Schulman_without_insertion_tolerance},
let $G^\dagger$ be its planar dual, let $X$ be an invariant bond percolation
process on $G$, let $X^\dagger$ be its dual, and let $k$ and
$k^\dagger$
be their (possibly random) number of infinite clusters. Then
Theorem \ref{thm:Newman-Schulman_without_insertion_tolerance} leaves
nine possibilities for the value of $(k, k^\dagger)\dvtx(0,0), (0,1),
(0, \infty), (1,0), (1,1), (1, \infty), (\infty, 0), (\infty, 1)$ and
$(\infty, \infty)$, but the following result rules out
four of them.
\begin{theorem}[(Theorem 3.1 in \cite{BS01})]
\label{thm:pairs_of_Ns}
With $G$, $G^\dagger$, $X$ and $X^\dagger$ as above, we have a.s. that
%
%
\begin{equation} \label{eq:pairs_of_Ns}
(k, k^\dagger) \in\{(0,1), (1,0), (1, \infty), (\infty, 1),
(\infty, \infty)\} .
\end{equation}
\end{theorem}

All five outcomes in (\ref{eq:pairs_of_Ns}) can actually happen.
The cases $(1, \infty)$ and $(\infty, 1)$ arise in the so-called
uniform spanning forest model (see Theorem \ref{thm:FSF_WSF_hyperbolic}
in the next section), but cannot arise in i.i.d. percolation
because of a local modification argument (Theorem 3.7 in \cite{BS01})
that can turn $(k, k^\dagger)=(1, \infty)$ into $(2, \infty)$,
contradicting Theorem \ref{thm:Newman-Schulman_without_insertion_tolerance}.
Benjamini and Schramm show that, in fact,
%
%
\begin{equation} \label{eq:pairs_in_iid_case}
\begin{tabular}{p{325pt}}
outcomes $(0,1)$, $(1,0)$ and $(\infty, \infty)$ are
exactly
those that happen for i.i.d. percolation,
\end{tabular}\hspace*{-28pt}
\end{equation}
and they prove the following remarkable result, establishing
Conjecture \ref{conj:nonuniqueness_phase} for the case of
planar hyperbolic graphs:\vspace*{-2pt}
\begin{theorem}[(Theorem 1.1 in \cite{BS01})] \label{thm:planar_nonuniqueness_phase}
Let $G=(V,E)$ be a planar nonamenable transitive graph with one end.
Then \mbox{$0 < p_{c,\mathrm{bond}}(G) < p_{u,\mathrm{bond}}(G) < 1$}. The same
is true for site percolation.\vspace*{-2pt}
\end{theorem}

This generalizes a result of Lalley \cite{Lal} who showed
\[
0 < p_{c,\mathrm{site}}(G) < p_{u,\mathrm{site}}(G) < 1
\]
for a more restrictive class of
graphs. Together with (\ref{eq:pairs_in_iid_case}),
this yields
\[
(k, k^\dagger) = \cases{
(0,1), &\quad for $p\in(0, p_{c,\mathrm{bond}}(G))$,\cr
(\infty, \infty), &\quad for $p \in(p_{c,\mathrm{bond}}(G), p_{u,\mathrm{bond}}(G))$,
\cr
(1,0), &\quad for $p \in(p_{u,\mathrm{bond}}(G), 1)$.}
\]
Concerning the behavior at the critical values,
Theorem \ref{thm:at_criticality} yields $(k, k^\dagger)=(0,1)$ for
$p=p_{c,\mathrm{bond}}(G)$, and by exchanging the roles of $G$ and $G^\dagger$
we see that $(k, k^\dagger)=(1,0)$ at $p=p_{u,\mathrm{bond}}(G)$.
As mentioned in Section \ref{sect:BS-beyond}, this
uniqueness of the infinite cluster already \textit{at} the uniqueness
critical value contrasts with the behavior obtained for certain
other nonamenable transitive graphs by Schonmann \cite{Scho99b}
and Peres \cite{Per}.

In fact, similar considerations for the hypothetical scenario
that Theorem~\ref{thm:planar_nonuniqueness_phase} fails show how
smoothly Theorem \ref{thm:planar_nonuniqueness_phase} follows from
(\ref{eq:pairs_in_iid_case}). We would then have $(k, k^\dagger)=(0,1)$
for $p<p_{c,\mathrm{bond}}(G)$ and $(k, k^\dagger)=(1,0)$
for $p>\break p_{c,\mathrm{bond}}(G)$. Hence $p_{c,\mathrm{bond}}(G^\dagger)=1-p_{c,\mathrm{bond}}(G)$, and
Theorem \ref{thm:at_criticality} applied to both $G$ and $G^\dagger$ yields
$(k, k^\dagger)=(0,0)$ at $p=p_{c,\mathrm{bond}}(G)$. This, however, contradicts
(\ref{eq:pairs_in_iid_case}).

Benjamini and Schramm \cite{BS01} go on to study properties of the infinite
clusters and their limit points on the boundary of the hyperbolic disk
(see also Lalley \cite{Lal,Lal01} for further results in this direction).
The final part of their paper \cite{BS01} concerns i.i.d. site
percolation on a certain random lattice in $\Hb^2$, namely
the Delaunay triangulation\vspace*{1pt} of a Voronoi tesselation for a~Poisson
process with intensity
$\lambda>0$ in~$\Hb^2$. The counterpart in $\R^d$ of such a~model
has also
been studied; see Bollob\'as and Riordan \cite{BRpaper,BR} for
a recent breakthrough in $\R^2$. However, the phase diagram becomes richer
and more interesting in the $\Hb^2$ case, not just because
of the nonuniqueness phase between $p_{c, \mathrm{site}}$ and $p_{u,\mathrm{site}}$,
but also because it becomes a true two-parameter family of models:
in $\R^d$, changing $\lambda$ is just a trivial rescaling of the model,
while in $\Hb^2$ there is no such scale invariance.

More recently, Tykesson \cite{Tyk} considered a different way of doing
percolation in $\Hb^2$ based on a Poisson process, namely to place a ball
of a fixed hyperbolic radius $R$ around each point, and consider
percolation properties of the region covered by the union of the balls.
(This is the so-called Boolean model of continuum percolation, which
has received a fair amount of attention in $\R^d$; see, e.g., Meester
and Roy \cite{MR}.) Equally natural is to consider percolation properties
of the vacant region (i.e., the complement of the covered region).
To consider percolation properties of both regions simultaneously is
analogous to working with the pair $(X, X^\dagger)$ in the discrete lattice
setting. The results in \cite{Tyk} turn out mostly
analogous to those of Benjamini and Schramm \cite{BS01} discussed
above for the discrete lattice setting.
However, in an even more recent paper by Benjamini, Jonasson,
Schramm and Tykesson \cite{BJST}, a phenomenon which is
particular to the continuum setting is revealed. Namely, are there
hyperbolic lines entirely contained in the vacant region, so that someone
living in $\Hb^2$ can actually ``see infinity'' in some
(necessarily random) directions? In $\R^d$ the answer to the analogous
question is no (this is related to Olbers' paradox in astronomy; see,
e.g., Harrison \cite{Harrison}) while in $\Hb^2$ the answer turns out
to be yes in certain parts of the parameter space.
The paper \cite{BJST}---which, sadly, became one of the
last by Oded Schramm---also contains results on various
refinements and variants of this question.

\section{Random spanning forests} \label{sect:forests}

So far, a lot has been said about percolation on nonamenable transitive
graphs under the general assumption of automorphism invariance, but
hardly anything about particular examples beyond the i.i.d. cases (other
than a few examples specifically designed to be counterexamples). But
as mentioned in Section \ref{sect:BLPS-gip}, there \textit{are} plenty
of important examples, and time has come to discuss one of them: the uniform
spanning forest. Later in this section, I will go
on to discuss its cousin, the \textit{minimal} spanning forest.

A spanning tree for a finite connected graph $G=(V,E)$ is a connected
subgraph containing all vertices but no cycles. A
\textit{uniform spanning tree}
for~$G$ is one chosen at random according to
uniform distribution on the set of possible spanning trees. Replacing
$G$ by, say, the $\Z^d$ lattice, the number of possible spanning
trees skyrockets to $\infty$, and it is no longer obvious how
to make sense of the uniform spanning tree. Pemantle \cite{Pem}
managed to make sense of it: for any infinite locally finite
connected graph $G=(V,E)$, let $\{G_1=(V_1, E_1), G_2=(V_2,
E_2),\vadjust{\eject}
\ldots\}$
be an increasing sequence of finite connected subgraphs of $G$, which
exhausts $G$ in the sense that each $e \in E$ and each $v \in V$ is in all
but at most finitely many of the $G_i$'s. Then, it turns out, the
uniform spanning tree measures for $G_1, G_2, \ldots$ converge weakly to
a probability measure $\mu_G$ on $\{0,1\}^E$ which is independent
of the exhaustion $\{G_1, G_2, \ldots\}$. In particular, $\mu_G$ is
$\Aut(G)$-invariant. Furthermore, it is concentrated on the event
that there are no open cycles and no finite open clusters, so that in
other words what we get is $\mu_G$-a.s. a forest, all of whose
trees are infinite. Naively, one might expect to get a single tree, but
this is not always the case. Pemantle showed
for the $\Z^d$ case that the number of trees is an a.s. constant
$N$ satisfying
%
%
\begin{equation} \label{eq:Pemantle_dichotomy}
N = \cases{
1, &\quad if $d\leq4$, \cr
\infty, &\quad otherwise.}
\end{equation}
This $d \leq4$ vs. $d>4$ dichotomy is related to the fact that
two independent SRW trajectories on the $\Z^d$ lattice intersect a.s.
if and only if $d \leq4$ (though this innocent-looking statement
hides the fact that Pemantle \cite{Pem} had to use deep results by
Lawler \cite{Law} on so-called loop-erased random walk; the proof
was simplified in the paper \cite{BLPSusf} to be discussed
below). The behavior for \mbox{$d>4$} suggests
\textit{uniform spanning forest} as a better term than uniform spanning
tree when $G$ is infinite.
The analysis of the uniform spanning forest~$\mu_G$ builds to a large extent
on the beautiful collection of identities
between uniform spanning trees, electrical networks and random walks which
has a~long and disperse history beginning with the
1847 paper by Kirchhoff \cite{Kir}. For instance, Rayleigh's
Monotonicity Principle for effective resistances in electrical networks
(see, e.g., Doyle and Snell \cite{DS}) underlies the stochastic monotonicity
properties of the sequence $\mu_{G_1}, \mu_{G_2}, \ldots$
that allows us to deduce the existence of the limiting measure $\mu_G$.

Some years after Pemantle's \cite{Pem} pioneering 1991 paper, and at about
the same time that BLPS \cite{BLPSgip} and BLPS \cite{BLPScritperc}
were written, the BLPS quartet started getting seriously interested in
uniform spanning forests; see Lyons~\cite{L08} for a more exact
statement about the timing and relation betweens the various BLPS projects.
This resulted in the magnificent paper \textit{Uniform spanning
forests} by
Benjamini, Lyons, Peres and Schramm \cite{BLPSusf}, which provides
a unified treatment and a number of simplications
of what was known on uniform spanning trees and forests,
together with a host of new and important results. The methods
involve, in addition to the aforementioned connections
to random walks and electrical networks, also mass-transport ideas,
Hilbert space projections, and Wilson's \cite{Wil} substantial
improvement on the Aldous--Broder algorithm \cite{Ald,Broder}
for generating uniform spanning trees. Here, I will mention only a couple
of the results from BLPS \cite{BLPSusf}, but see Lyons \cite{L98} for
a gentle introduction to the same topic (note that despite the inverted
publication dates, the 1998 paper \cite{L98} surveys much of the
original work in the 2001 paper \cite{BLPSusf}).

It turns out that there is another, equally natural, way to obtain
a uniform spanning forest for an infinite graph $G=(V,E)$ via
the exhaustion $\{G_1, G_2, \ldots\}$ considered
above, namely if for each $i$ we consider the uniform spanning tree
not on $G_i$ but on the modified graph where an extra vertex $w$ is
introduced, together with edges $\langle v,w \rangle$ for all
$v\in V_i \setminus V_{i-1}$. We think of this as a ``wired'' version
of $G_i$, and the limiting measure on $\{0,1\}^E$, which is denoted
$\WSF_G$ and which exists for similar reasons as for $\mu_G$, should
be thought of as the \textit{wired uniform spanning forest} for~$G$. For
consistency of terminology, we call $\mu_G$ the \textit{free uniform spanning
forest}, and switch to denoting it $\FSF_G$. (The free and wired
uniform spanning forests are analogous to the free and
wired limiting measures for the random-cluster model; cf. \cite{Gr06}.)
The wired measure $\WSF_{\Z^d}$ appeared implicitly in Pemantle \cite{Pem}
together with the result that $\WSF_{\Z^d}=\FSF_{\Z^d}$; this was made
explicit in H\"aggstr\"om~\cite{H95}, and BLPS~\cite{BLPSusf} noted
that this extends to $\WSF_G = \FSF_G$ whenever~$G$ is transitive
and amenable. On the other hand, $\WSF_{\T_d}\neq\FSF_{\T_d}$ for
the $(d+1)$-regular tree $\T_d$ when $d \geq2$, so a first guess
might be that for transitive graphs, $\WSF_G = \FSF_G$ is equivalent
to amenability. This turns out not to be true, however,
and BLPS \cite{BLPSusf} offer instead the following characterization
(which does not require $G$ to be transitive or even quasi-transitive).
Recall that for a graph $G=(V,E)$ a function $f\dvtx V \rightarrow\R$
is called \textit{harmonic} if for any $v \in V$ we have that $f(v)$
equals the average of $f(w)$ among all its neighbors $w$, and that
$f$ is called \textit{Dirichlet} if
$\sum_{\langle u,v \rangle\in E}(f(u)-f(v))^2<\infty$.\looseness=1
\begin{theorem}[(Theorem 7.3 in \cite{BLPSusf})] \label{thm:FSF=WSF}
For any graph $G$, we have $\WSF_G = \FSF_G$ if and only if $G$
admits no
nonconstant harmonic Dirichlet functions.
\end{theorem}

As an example of a nonamenable transitive graph for which $\WSF_G =
\FSF_G$
holds, we may take the Grimmett--Newman example discussed in
Section \ref{sect:BS-beyond}; this follows from
Theorem \ref{thm:FSF=WSF} in
combination with the result of Thomassen \cite{Thom,Soardi}
that the Cartesian product of any two infinite graphs has no nonconstant
harmonic Dirichlet functions.

A nice class of graphs where $\WSF_G \neq\FSF_G$ are the
planar hyperbolic lattices considered in the previous section:
\begin{theorem}[(Theorems 12.2 and 12.7 in BLPS \cite{BLPSusf})] \label
{thm:FSF_WSF_hyperbolic}
For any planar nonamenable transitive graph $G=(V,E)$ with one end, we have
that $\WSF_G \neq\FSF_G$. In particular, they differ in the number
of infinite clusters: $\FSF_G$ produces a unique infinite
cluster a.s., while $\WSF_G$ produces infinitely many a.s.
If $G^\dagger=(V^\dagger, E^\dagger)$
is the planar dual of $G$, and $X \in\{0,1\}^E$ is given
by $\WSF_G$, then the dual percolation $X^\dagger\in\{0,1\}
^{E^\dagger}$
defined as in Section \ref{sect:planar}, has
distribution $\FSF_{G^\dagger}$.
\end{theorem}

There is a lot more to quote from BLPS \cite{BLPSusf} concerning
uniform spanning forests, but I will instead cite a result
from a later paper by Benjamini, Kesten, Peres and Schramm \cite{BKPS}.
Strictly speaking the result falls outside the scope of the present
survey as I defined it in Section \ref{sect:intro},
because it concerns~$\Z^d$, but it is so
supremely beautiful that I find this inconsistency of mine
to be motivated. For a bond percolation process $X\in\{0,1\}^E$
on a graph $G=(V,E)$, and two vertices $u,v\in V$,
define the random variable $D_{\mathrm{closed}}(u,v)$ as the minimal number
of closed edges that a path in $G$ from $u$ to $v$ needs to traverse,
and $D_{\mathrm{closed}}^{\max}$ as the supremum of $D_{\mathrm{closed}}(u,v)$ over
all choices of $u,v\in V$. For percolation processes such as
$\WSF_G$ and $\FSF_G$ where a.s. all vertices belong to infinite clusters,
$D_{\mathrm{closed}}^{\max}=0$ means precisely uniqueness of the infinite cluster,
while $D_{\mathrm{closed}}^{\max}=1$ means that uniqueness fails but that any
pair of infinite clusters come within unit distance from each other somewhere
in $G$. The dichotomy by Pemantle quoted in
(\ref{eq:Pemantle_dichotomy}) says that for uniform spanning forests
on $\Z^d$, we have $D_{\mathrm{closed}}^{\max}=0$ a.s. for $d\leq4$ but
$D_{\mathrm{closed}}^{\max} \geq1$ a.s. when $d >4$. Benjamini, Kesten, Peres and
Schramm \cite{BKPS} proved the following refinement.
\begin{theorem} \label{thm:BKPS}
For the uniform spanning forest measure $\WSF_{\Z^d}$ (or equivalently
$\FSF_{\Z^d}$), we have, a.s.,
\[
D_{\mathrm{closed}}^{\max} = \cases{
0, &\quad if $d \in\{1,2,3,4\}$, \cr
1, &\quad if $d \in\{5,6,7,8\}$, \cr
2, &\quad if $d \in\{9,10,11,12\}$, \cr
3, &\quad if $d \in\{13,14,15,16\}$, \cr
&\quad \vdots
}
\]
\end{theorem}

Next, let us return briefly to the case
where $G=(V,E)$ is a finite graph. Besides
the uniform spanning tree, there is another, much-studied
and equally natural,
way of picking a random spanning tree for $G$: attach i.i.d. weights
$\{U(e)\}_{e \in E}$ to the edges of $G$,
and pick the spanning tree for $G$ that minimizes the sum of the
edge weights; this is the so-called \textit{minimal spanning tree} for~$G$. The
marginal distribution of the $U(e)$'s has no effect on the
distribution of the tree as long as it is free from atoms, but it turns out
to facilitate the analysis to take it to be uniform on $[0,1]$.
It is easy to see that an edge $e \in E$ is included of the minimal
spanning tree if and only if
$U(e)< \max_{e \in C}U(e')$ for all cycles $C$ containing $e$.

If instead $G$ is infinite, this
characterization can be taken as the definition, and the resulting
random subgraph of $G$ is called the \textit{free minimal spanning forest},
the corresponding probability measure on $\{0,1\}^E$ being denoted
$\FMSF_G$. An alternative extension to infinite $G$ is to include
$e\in E$ if and only if
$U(e)< \max_{e \in C}U(e')$ for all \textit{generalized}
cycles $C$ containing $e$, where by a generalized cycle we mean
a cycle or a bi-infinite self-avoiding path.
This gives rise to the so-called \textit{wired minimal spanning forest},\vadjust{\eject}
and a corresponding probability measure $\WMSF_G$ on $\{0,1\}^E$. The
study of $\FMSF_G$ and $\WMSF_G$ parallels the study of $\FSF_G$ and
$\WSF_G$
in many ways, and I wish to draw the reader's attention to the paper
\textit{Minimal spanning forests} by Lyons, Peres and Schramm \cite{LPS},
where, in the words of the authors, ``a key theme [is] to describe
striking analogies, and important differences, between uniform and
minimal spanning forests.''

As for uniform spanning forests, one of the central issues is to determine
when $\FMSF_G = \WMSF_G$. And as for uniform spanning forests,
equality turns
out to hold on $\Z^d$ and more generally on amenable transitive graphs.
But moving on to the nonamenable case, the statements $\FMSF_G = \WMSF_G$
and $\FSF_G=\WSF_G$ are no longer equivalent.
\begin{theorem}[(Proposition 3.6 in \cite{LPS})] \label{thm:MSF_uniqueness}
On any $G$, we have $\FMSF_G = \WMSF_G$ if and only if
for Lebesgue-a.e. $p\in[0,1]$ it is the case that i.i.d. bond
percolation on $G$ with parameter $p$
produces a.s. at most one infinite cluster.
\end{theorem}

Hence we should expect to have
$\FMSF_G \neq\WMSF_G$ for nonamenable quasi-transitive $G$; this
is equivalent to Conjecture \ref{conj:nonuniqueness_phase}. One example
where we do know that $\FMSF_G \neq\WMSF_G$ while $\FSF_G=\WSF_G$ is
when $G$ is the Grimmett--Newman graph $\T_d \times\Z$; cf.
Theorem \ref{thm:FSF=WSF} and the comment following it.

Another
consequence of having a nonuniqueness phase for i.i.d. bond percolation
on $G$ is that uniqueness of the infinite cluster fails for
$\WMSF_G$ (this is Corollary 3.7 in \cite{BLS}). Determining the number of
infinite clusters in $\WMSF_G$ (or in $\FMSF_G$)
more generally is of course a central issue.
On~$\Z^d$ for $d \geq2$, the only case that has been settled is $d=2$,
where Alexander~\cite{Alex} showed that the number of infinite clusters
is $1$. For higher dimensions, a~dichotomy like Pemantle's
(\ref{eq:Pemantle_dichotomy}) can be expected; Newman and Stein \cite{NSt}
conjectured that the switch from uniqueness to infinitely
many infinite clusters should happen at $d=8$ or $9$.

A general difference---but at the same time a parallel---between
the uniform spanning forest and the minimal
spanning forest is that while the former has intimate connections
to SRW on $G$, the latter seems equally intimately connected to
i.i.d. bond percolation on $G$, such as in Theorem \ref{thm:MSF_uniqueness}.
[The key device for exposing the connections between the minimal
spanning forest and i.i.d. percolation is the
coupling in which we use the $U(e)$'s underlying
$\FMSF_G$ and $\WMSF_G$ also for generating i.i.d. bond percolation
for any $p$: an edge $e$ is declared open on level $p$ iff $U(e)<p$.]
One striking instance of this parallel is the following. Morris \cite{Mor}
showed for any $G$ that a.s. any component $C$ of the wired uniform
spanning forest has the property that SRW on it is recurrent, while
Lyons, Peres and Schramm \cite{LPS} showed for any $G$ that
a.s. any component $C$ of the wired minimal
spanning forest has $p_{c,\mathrm{bond}}(C)=1$. Both results are sharp in
the sense that in neither of them can the set of possible
$C$'s be narrowed further without making the statement false.

\section{Postscript} \label{sect:final_word}

The reader may have noticed that most of the work surveyed
here is from the mid-to-late 1990s, and ask why this is.
Is it because Oded got bored with percolation beyond $\Z^d$? No, it is not,
and no, he did not. The main reason is Oded's discovery in \cite{S00}
of SLE
(Schramm--Loewner evolution, or stochastic Loewner evolution as Oded
himself preferred to call it). This led him to the far-reaching
and more urgent project, beginning with a now-famous series of
papers with Greg Lawler and Wendelin Werner \cite
{LSW1,LSW2,LSW3,LSW4}, of
understanding SLE and how it arises as a scaling limit of
various critical models in two dimensions. But that is a different story.

\section*{Acknowledgment}
I am grateful to Russ Lyons and Johan Tykesson for
helpful comments on an earlier draft.

%

%
\printaddresses

\end{document}